\documentclass{article}
\usepackage{amsmath, latexsym, amsfonts, amssymb, amsthm, amscd}

\newtheorem {lemme} {Lemma} [section]
\newtheorem {theoreme} {Theorem} [section]
\newtheorem {proposition} {Proposition} [section]
\newtheorem {corollaire} {Corollary} [section]
\newtheorem {remarque} {Remark} [section]

\newcommand{\tr}{{\rm tr}}
\newcommand{\Tr}{{\rm Tr}}
\newcommand{\E}{\mathbb {E}}
\newcommand{\N}{\mathbb {N}}
\newcommand{\R}{\mathbb {R}}
\newcommand{\C}{\mathbb {C}}

\newcommand{\1}{1\!\!{\sf I}}

\numberwithin{equation}{section}

\newcommand{\vers}{\mathop{\longrightarrow }}
\title{ Free convolution with a semicircular distribution and 
eigenvalues of spiked deformations of Wigner matrices
\footnote{This work was partially supported by the {\emph Agence Nationale de la
Recherche} grant ANR-08-BLAN-0311-03.}}

\author{M. Capitaine\thanks{CNRS, Institut de Math\'ematiques de Toulouse, 
Equipe de Statistique et Probabilit\'es,  F-31062 Toulouse Cedex 09. 
E-mail: mireille.capitaine@math.univ-toulouse.fr },
\hspace{.1cm}
C. Donati-Martin\thanks{Universit\'e de Versailles-St Quentin, 
Laboratoire de Math\'ematiques,
45 avenue de Etats Unis, F-78035 Versailles Cedex. E-mail: catherine.donati-martin@uvsq.fr }
\hspace{.1cm},
D. F\'eral\thanks{Institut de Math\'ematiques de Bordeaux, 
Universit\'e Bordeaux 1, 351 Cours de la Lib\'eration, F-33405 Talence Cedex. 
E-mail: delphine.feral@math.u-bordeaux1.fr } 
\hspace{.1cm} 
and M. F\'evrier\thanks{Institut de Math\'ematiques de Toulouse, 
Equipe de Statistique et Probabilit\'es,  F-31062 Toulouse Cedex 09. 
E-mail: fevrier@math.univ-toulouse.fr }}

\date{}
\begin{document}
\maketitle
\begin{abstract}
We investigate the asymptotic behavior of the eigenvalues of spiked perturbations of Wigner matrices defined by 
$M_N=\frac{1}{\sqrt{N}}W_N+A_N$, where 
$W_N$ is a $N\times N$ Wigner Hermitian matrix whose entries have
a distribution $\mu $ which is symmetric and satisfies a Poincar\'e inequality and 
$A_N$ is a deterministic Hermitian matrix whose spectral measure 
converges to some probability measure $\nu $ with compact support.
We assume that $A_N$ has a fixed number of fixed eigenvalues (spikes) 
outside the support of $\nu $ whereas the distance between the other eigenvalues 
and the support of $\nu $ uniformly goes to zero as $N$ goes to infinity. 
We establish that only a particular subset of the spikes will generate 
some eigenvalues of $M_N$ which will converge to some limiting points 
outside the support of the limiting spectral measure.
This phenomenon can be fully described in terms of free probability 
involving the subordination function related 
to the free additive convolution of $\nu $ by a semicircular distribution.
Note that only finite rank perturbations had been considered up to now 
(even in the deformed GUE case). 
\end{abstract}



\noindent Submitted to EJP on  November 8, 2010, final version accepted on  August 9, 2011.
\section{Introduction}\label{intro}

In the fifties, in order to describe the energy levels of a complex nuclei system 
by the eigenvalues of large Hermitian matrices, 
E. Wigner introduced the so-called Wigner $N\times N$ matrix $W_N$. 
According to Wigner's work \cite{Wigner55}, \cite{Wigner58} 
and further results of different authors (see \cite{Bai99} for a review), 
provided the common distribution $\mu $ of the entries is centered
with variance $\sigma ^2$, the large $N$-limiting spectral distribution 
of the rescaled complex Wigner matrix $X_N=\frac{1}{\sqrt{N}} W_N$ is 
the semicircle distribution $\mu _{\sigma }$ whose density is given by
\begin{equation}\label{scl}
\frac{d\mu _{\sigma }}{dx}(x)= \frac{1}{2 \pi \sigma ^2} \sqrt{4\sigma ^2- x^2} 
\, 1 \hspace{-.20cm}1_{[-2\sigma , 2\sigma ]}(x).
\end{equation}
Moreover, if the fourth moment of the measure $\mu $ is finite, 
the largest (resp. smallest) eigenvalue of $X_N$ 
converges almost surely towards the right (resp. left) endpoint $2\sigma $ 
(resp. $-2\sigma $) of the semicircular support 
(cf. \cite{BaiYin88} or Theorem 2.12 in \cite{Bai99}). 

Now, how does the spectrum behave under a deterministic Hermitian perturbation $A_N$? 
The set of possible spectra for $M_N=X_N + A_N$ depends in a complicated way 
on the spectra of $X_N$ and $A_N$ (see \cite{Fulton98}). 
Nevertheless, when $N$ becomes large, free probability provides us a good understanding 
of the global behavior of the spectrum of $M_N$. 
Indeed, if the spectral measure of $A_N$ weakly converges to some probability measure $\nu$ 
and $\Vert A_N\Vert $ is uniformly bounded in $N$, 
the spectral distribution of $M_N$ weakly converges to the free convolution $\mu _\sigma \boxplus \nu $ 
almost surely and in expectation 
(cf \cite{AGZ09}, \cite{MinSpe10} and \cite{Voiculescu91}, \cite{Dykema93} for pioneering works). 
We refer the reader to \cite{VDN92} for an introduction to free probability theory. 
Note that when $A_N$ is of finite rank, 
the spectral distribution of $M_N$ still converges to the semicircular distribution 
($\nu\equiv \delta_0$ and $\mu _\sigma \boxplus \nu =\mu _\sigma $).\\
In \cite{Peche06}, S. P\'ech\'e investigated the deformed GUE model $M_N^G=W_N^G/\sqrt{N}+A_N$, 
where $W_N^G$ is a GUE matrix, that is a Wigner matrix 
associated to a centered Gaussian measure with variance $\sigma ^2$ 
and $A_N$ is a deterministic perturbation of finite rank with fixed eigenvalues. 
This model is the additive analogue of the Wishart matrices 
with spiked covariance matrix previously considered by 
J. Baik, G. Ben Arous and S. P\'ech\'e \cite{BBP05} 
who exhibited a striking phase transition phenomenon 
for the fluctuations of the largest eigenvalue according to the values of the spikes. 
S. P\'ech\'e pointed out an analogous phase transition phenomenon 
for the fluctuations of the largest eigenvalue of $M_N^G$ 
with respect to the largest eigenvalue $\theta $ of $A_N$ \cite{Peche06}. 
These investigations imply that, if $\theta $ is far enough from zero ($\theta > \sigma $), 
then the largest eigenvalue of $M_N^G$ jumps above the support $[-2\sigma , 2\sigma ]$ 
of the limiting spectral measure and converges (in probability) 
towards $\rho _\theta =\theta + \frac{\sigma ^2}{\theta }$. 
Note that Z. F\"uredi and J. Koml$\acute{\text{o}}$s already exhibited such a phenomenon
in \cite{FurKom81} dealing with non-centered symmetric matrices.\\
In \cite{FP07}, D. F\'eral and S. P\'ech\'e proved that the results of \cite{Peche06} 
still hold for a non-necessarily Gaussian Wigner Hermitian matrix $W_N$ 
with sub-Gaussian moments and in the particular case of a rank one perturbation matrix $A_N$ 
whose entries are all $\frac{\theta }{N}$ for some real number $\theta $. 
In \cite{CDF09}, we considered a deterministic Hermitian matrix $A_N$ 
of arbitrary fixed finite rank $r$ 
and built from a family of $J$ fixed non-null real numbers 
$\theta _1 > \cdots > \theta _J$ independent of $N$ 
and such that each $\theta _j$ is an eigenvalue of $A_N$ of fixed multiplicity $k_j$ 
(with $\sum_{j=1}^J k_j=r$). 
In the following, the $\theta_j$'s are referred as the spikes of $A_N$. 
We dealt with general Wigner matrices associated to some symmetric measure satisfying a Poincar\'e inequality. 
We proved that eigenvalues of $A_N$ with absolute value strictly greater than $\sigma $
generate some eigenvalues of $M_N$ which converge 
to some limiting points outside the support of $\mu _\sigma $. 
To be more precise, we need to introduce further notations. 
Given an arbitrary Hermitian matrix $B$ of size $N$, 
we denote by $\lambda _1(B) \geq \cdots \geq \lambda _N(B)$
its $N$ ordered eigenvalues. For each spike $\theta _j$, 
we denote by $n_{j-1}+1, \ldots, n_{j-1}+k_j$ the descending ranks of $\theta _j$ 
among the eigenvalues of $A_N$ (multiplicities of eigenvalues are counted) 
with the convention that $k_1+ \cdots +k_{j-1}=0$ for $j=1$. One has that 
$$n_{j-1}=k_1+ \cdots +k_{j-1} \quad \text{if $\theta _j > 0$} \quad \text{and} 
\quad n_{j-1}=N-r+k_1+ \cdots +k_{j-1} \quad \text{if $\theta _j < 0$}.$$
Letting $J_{+\sigma }$ (resp. $J_{-\sigma }$) be the number of $j$'s such 
that $\theta _j > \sigma $ (resp. $\theta _j < -\sigma $), 
we established in \cite{CDF09} that, when $N$ goes to infinity, 
\begin{itemize}
\item[{a)}] for all $j$ such that $1 \leq j \leq J_{+\sigma }$ 
(resp. $j \geq J-J_{-\sigma }+1$), the $k_j$ eigenvalues 
$(\lambda _{n_{j-1}+i}(M_N), 1 \leq i \leq k_j)$ converge almost surely 
to $\rho _{\theta _j}=\theta _j+\frac{\sigma ^2}{\theta _j}$ 
which is $> 2\sigma $ (resp. $< -2\sigma $). 
\item[{b)}] $\lambda _{k_1+ \cdots + k_{J_{+\sigma }}+1}(M_N)\overset{a.s.}{\longrightarrow }2\sigma $ 
\text{and} $\lambda _{N-(k_J+ \cdots + k_{J-J_{-\sigma }+1})}(M_N)\overset{a.s.}{\longrightarrow }-2\sigma $.
\end{itemize}

Actually, this phenomenon may be described in terms of free probability 
involving the subordination function related to the free convolution 
of $\nu =\delta _0$ by a semicircular distribution. Let us present it briefly.
For a probability measure $\tau $ on $\R$, 
let us denote by $g_\tau$ its Stieltjes transform, defined for $z\in \C\setminus \R$ by
$$g_\tau (z) = \int_\R \frac{d\tau(x)}{z-x}.$$
Let $\nu $ and $\tau $ be two probability measures on $\R$. 
It is proved in \cite{Biane98} Theorem 3.1 that there exists an analytic map 
$F : \C^+ \rightarrow \C^+$, called subordination function, 
such that $$\forall z \in \C^+ , g_{\tau \boxplus \nu }(z)= g_{\nu }(F(z)),$$
where $\C^+$ denotes the set of complex numbers $z$ such that $\Im z> 0$. 
When $\tau = \mu _\sigma $, let us denote by $F_{\sigma ,\nu }$ the corresponding subordination function.
When $\nu = \delta _0$ and $\tau =\mu _\sigma $, 
the subordination function is given by $F_{\sigma ,\delta _0} = 1/g_{\mu _\sigma }$. 
According to Lemma 4.4 in \cite{CDF09}, one may notice that the complement of the support 
of $\mu _\sigma \boxplus \delta _0 (= \mu _\sigma )$ can be described as: 
$$\R\setminus [-2\sigma , 2\sigma ] = \{x, \exists u \in \R^*, \vert u\vert > \sigma 
~~\mbox{such that~~} x = H_{\sigma , \delta _0}(u) \},$$
where $H_{\sigma ,\delta _0}(z) = z + \frac{\sigma^2}{z}$ 
is the inverse function of the subordination function 
$F_{\sigma , \delta _0}$ on $\R\setminus [-2\sigma , 2\sigma ]$. 
Now, the characterization of the spikes of $A_N$ 
that generate jumps of eigenvalues of $M_N$ i.e. $\vert \theta _j\vert > \sigma $ 
is obviously equivalent to the following 
$$\theta_j \in \R\setminus {\rm supp}(\delta _0) (=\R^*) 
\quad \text{and} \quad H_{\sigma , \delta _0}'(\theta_j) > 0.$$
Moreover the relationship between a spike $\theta_j$ of $A_N$ 
such that $\vert \theta_j\vert > \sigma $ and the limiting point 
$\rho _{\theta _j}$ of the corresponding eigenvalues of $M_N$ (which is then outside $[-2\sigma ; 2\sigma ]$) 
is actually described by the inverse function of the subordination function as: 
$$\rho _{\theta _j}=H_{\sigma , \delta _0}(\theta _j).$$

\noindent 
Actually this very interpretation in terms of subordination function 
of the characterization of the spikes of $A_N$ that generate jumps of eigenvalues of $M_N$ 
as well as the values of the jumps provides the intuition to imagine 
the generalization of the phenomenon dealing with non-finite rank perturbations 
just by replacing $\delta _0$ by the limiting spectral distribution $\nu $ of $A_N$ in the previous lines. 
Up to now, no result has been established for non-finite rank additive spiked perturbation. 
Moreover, this paper shows up that free probability can also shed light 
on the asymptotic behavior of the eigenvalues of the deformed Wigner model 
and strengthens the fact that free probability theory and random matrix theory are closely related.

More precisely, in this paper, we consider the 
following general deformed Wigner models $M_N=X_N+ A_N$ such that: 
\begin{itemize}
\item $X_N=\frac{1}{\sqrt{N}} W_N $ where $W_N$ is a $N \times N$ Wigner Hermitian matrix associated to 
a distribution $\mu $ of variance $\sigma ^2$ and mean zero: \\
$(W_N)_{ii}$, $\sqrt{2} \Re ((W_N)_{ij})_{i < j}$, $\sqrt{2}
\Im ((W_N)_{ij})_{i<j}$
are  i.i.d., with distribution $\mu $ which is symmetric 
and satisfies a Poincar\'e inequality (the definition of 
such an inequality is recalled in the Appendix). 
\item $A_N$ is a deterministic Hermitian matrix whose eigenvalues $\gamma ^{(N)}_i$, 
denoted for simplicity by $\gamma _i$, 
are such that the spectral measure $\mu _{A_N} := \frac{1}{N} \sum_{i=1}^N \delta _{\gamma _i}$ 
converges to some probability measure $\nu $ with compact support. 
We assume that there exists a fixed integer $r\geq 0$ (independent from $N$) 
such that $A_N$ has $N-r$ eigenvalues $\beta _j(N)$ satisfying 
$$\max _{1\leq j\leq N-r} {\rm dist}(\beta _j(N),{\rm supp}(\nu ))\vers _{N \rightarrow \infty } 0,$$
where ${\rm supp}(\nu )$ denotes the support of $\nu $. 
We also assume that there are $J$ fixed real numbers $\theta _1 > \ldots > \theta _J$ 
independent of $N$ which are outside the support of $\nu $ and such that each $\theta _j$ 
is an eigenvalue of $A_N$ with a fixed multiplicity $k_j$ (with $\sum_{j=1}^J k_j=r$). 
The $\theta_j$'s will be called the spikes or the spiked eigenvalues of $A_N$. 
\end{itemize} 

\noindent According to \cite{AGZ09}, the spectral distribution of $M_N$ 
weakly converges to the free convolution $\mu _\sigma \boxplus \nu $ almost surely 
(cf. Remark \ref{SpectralASCV} below). 
It turns out that the spikes of $A_N$ that will generate jumps of eigenvalues of $M_N$ 
will be the $\theta _j$'s such that $H_{\sigma , \nu }'(\theta_j) > 0$ where 
$H_{\sigma , \nu }(z)=z+\sigma ^2g_\nu (z)$ and the corresponding limiting points 
outside the support of $\mu _\sigma \boxplus \nu $ will be given by
$$\rho _{\theta _j}=H_{\sigma , \nu }(\theta_j).$$
It is worth noticing that the set $\{u \in \R\setminus {\rm supp}(\nu ), \, H_{\sigma , \nu }'(u) > 0\}$ 
is actually the complement of the closure of the open set 
$$U_{\sigma , \nu }:= \left\{ u \in \R, \int_{\R} \frac{d\nu (x)}{(u-x)^2} > \frac{1}{\sigma ^2} \right\}$$
introduced by P. Biane in \cite{Biane97b} to describe the support of the free additive convolution 
of a probability measure $\nu $ on $\R$ by a semicircular distribution. 
Note that the deep study by P. Biane of the free convolution by a semicircular distribution 
will be of fundamental use in our approach. 
In Theorem \ref{ThmASCV}, which is the main result of the paper, 
we present a complete description of the convergence of the eigenvalues of $M_N$ 
depending on the location of the $\theta _j$'s with respect to $\overline{U_{\sigma , \nu }}$ 
and to the connected components of the support of $\nu $.
~

Our approach also allows us to study the ``non-spiked" deformed Wigner models i.e. such that $r=0$. 
Up to now, the results which can be found in the literature for such a situation 
concern the so-called Gaussian matrix models with external source 
where the underlying Wigner matrix is from the GUE. 
Many works on these models deal with the local behavior of the eigenvalues of $M_N$ 
(see for instance \cite{BK04}, \cite{ABK05} and \cite{BK07} for details). 
Moreover, the recent results of \cite{Male10} 
(which investigate several matrices in a free probability context) 
imply that the operator norm (i.e. the largest singular value) 
of some non-spiked deformed GUE $M_N^G=W_N^G/N + A_N$ converges almost surely 
to the $L^\infty $-norm of a $(\mu _\sigma \boxplus \nu )$-distributed random variable. 
Here, we readily deduce (cf. Proposition \ref{ThmASCVNonSpike} below) 
from our results the almost sure convergence of the extremal eigenvalues 
of general non-spiked deformed Wigner models to the corresponding endpoints 
of the compact support of the free convolution $\mu _\sigma \boxplus \nu $.

The asymptotic behavior of the eigenvalues of the deformed Wigner model $M_N$
actually comes from two phenomena involving free convolution:
\begin{enumerate}
\item {the inclusion of the spectrum of $M_N$ in an $\epsilon $-neighborhood 
of the support of $\mu _\sigma \boxplus \mu _{A_N}$, for all large $N$ almost surely;}
\item {an exact separation phenomenon between the spectrum of $M_N$ and the spectrum of $A_N$, 
involving the subordination function $F_{\sigma , \nu }$ of $\mu _\sigma \boxplus \nu $ 
(i.e. to a gap in the spectrum of $M_N$, it corresponds through $F_{\sigma , \nu }$ 
a gap in the spectrum of $A_N$ which splits the spectrum of $A_N$ exactly as that of $M_N$).}
\end{enumerate}
The key idea to prove the first point is to obtain a precise estimate of order $\frac{1}{N}$ 
of the difference between the respective Stieltjes transforms of the mean spectral measure 
of the deformed model and of $\mu _\sigma \boxplus \mu _{A_N}$. 
To get such an estimate, we prove an ``approximative subordination equation" 
satisfied by the Stieltjes transform of the deformed model. 
Note that, even if the ideas and tools are very close to those developed in \cite{CDF09}, 
the proof in \cite{CDF09} does not use the above analysis from free probability 
whereas this very analysis allows us to extend the results of \cite{CDF09} 
to non-finite rank deformations. 
In particular, we didn't consider in \cite{CDF09} $\mu _\sigma \boxplus \mu _{A_N}$ 
whose support actually makes the asymptotic values of the eigenvalues 
that will be outside the limiting support of the spectral measure of $M_N$ appear. 

Note that phenomena 1. and 2. are actually the additive analogues 
of those described in \cite{BaiSil98}, \cite{BaiSil99} in the framework of spiked population models, 
even if the authors do not refer to free probability. 
In \cite{BaikSil06}, the authors use the results of \cite{BaiSil98}, \cite{BaiSil99} 
to establish the almost sure convergence of the eigenvalues 
generated by the spikes in a spiked population model 
where all but finitely many eigenvalues of the covariance matrix are equal to one. 
Thus, they generalize the pioneering result of \cite{BBP05} in the Gaussian setting. 
Recently, \cite{RaoSil09}, \cite{BaiYao08b} extended this theory to 
a generalized spiked population model where the base population covariance matrix is arbitrary. 
Our results are exactly the additive analogues of theirs. 
It is worth noticing that one may check that these results on spiked population models 
could also be fully described in terms of free probability involving the subordination function 
related to the free multiplicative convolution of $\nu $ by a Marchenko-Pastur distribution. \\
Moreover, the results of F. Benaych-Georges and R. R. Nadakuditi in \cite{BGRao09} 
about the convergence of the extremal eigenvalues of a matrix $X_N +A_N$, 
$A_N$ being a finite rank perturbation whereas $X_N$ is a unitarily invariant matrix 
with some compactly supported limiting spectral distribution $\mu $, 
could be rewritten in terms of the subordination function related to 
the free additive convolution of $\delta _0$ by  $\mu $. 
Hence, we think that subordination property in free probability 
definitely sheds light on spiked deformed models.

Finally, one can expect that our results hold true in a more general setting 
than the one considered here, namely only requires the existence 
of a finite fourth moment on the measure $\mu $ of the Wigner entries. 
Nevertheless, the assumption that $\mu $ satisfies a Poincar\'e inequality 
is fundamental in our approach since we need several variance estimates. 

The paper is organized as follows. In Section \ref{freeconv}, 
we first recall some results on free additive convolution 
and subordination property as well as the description by P. Biane 
of the support of the free convolution of some probability measure $\nu $ 
by a semicircular distribution. 
We then deduce a characterization of this support via the subordination function 
when $\nu $ is compactly supported and we exhibit relationships between 
the steps of the distribution functions of $\nu $ and $\mu _\sigma \boxplus \nu $. 
In Section \ref{sec2}, we establish an approximative subordination equation 
for the Stieltjes transform $g_N$ of the mean spectral distribution of the deformed model $M_N$ 
and explain in Section \ref{difference} how to deduce an estimation up to the order $\frac{1}{N^2}$ 
of the difference between $g_N$ and the Stieltjes transform of $\mu _\sigma \boxplus \mu _{A_N}$ 
when $N$ goes to infinity. 
In Section \ref{inclusionN}, we show how to deduce the almost sure inclusion 
of the spectrum of $M_N$ in a neighborhood of the support of 
$\mu _\sigma \boxplus \mu _{A_N}$ for all large $N$; 
we use the ideas (based on inverse Stieltjes tranform) of \cite{HaaThor05} and \cite{Schultz05} 
in the non-deformed Gaussian complex, real or symplectic Wigner setting; 
nevertheless, since $\mu _\sigma \boxplus \mu _{A_N}$ depends on $N$, we need here 
to apply the inverse Stieltjes tranform to functions depending on $N$ 
and we therefore give the details of the proof to convince the reader 
that the approach developped by \cite{HaaThor05} and \cite{Schultz05} still holds. 
In Section \ref{supportN}, we show how the support of $\mu _\sigma \boxplus \mu _{A_N}$ 
makes the asymptotic values of the eigenvalues that will be outside the support 
of the limiting spectral measure appear since we prove that, for any $\epsilon > 0$, 
${\rm supp}(\mu _\sigma \boxplus \mu _{A_N})$ is included in an 
$\epsilon$-neighborhood of ${\rm supp}(\mu _\sigma \boxplus \nu ) \bigcup 
\left\{ \rho _{\theta _j}, \, \theta_j \mbox{~such that~}  \, H_{\sigma , \nu }'(\theta _j) > 0 \right\}$, 
when $N$ is large enough. 
Section \ref{exact} is devoted to the proof of the exact separation phenomenon 
between the spectrum of $M_N$ and the spectrum of $A_N$, 
involving the subordination function $F_{\sigma , \nu }$. 
In the last section, we show how to deduce our main result (Theorem \ref{ThmASCV}) 
about the convergence of the eigenvalues of the deformed model $M_N$. 
Finally we present in an Appendix the proofs of some technical estimates on variances used throughout the paper. 

Throughout this paper, we will use the following notations.
\begin{itemize}
\item[-] For a probability measure $\tau $ on $\R$, 
we denote by $g_\tau $ its Stieltjes transform defined for $z \in \C\setminus \R$ by 
$$g_\tau (z) = \int_\R \frac{d\tau (x)}{z-x}.$$
\item[-] $G_N$ denotes the resolvent of $M_N$ and 
$g_N$ the mean of the Stieltjes transform of the spectral measure of $M_N$, that is 
$$g_N(z) = \E(\tr_N G_N(z)), \, z \in \C\setminus \R,$$
where $\tr_N$ is the normalized trace: $\tr_N =\frac{1}{N} \Tr$.\\
We recall some useful properties of the resolvent (see \cite{KKP96}, \cite{CD07}).
\begin{lemme} \label{lem0}
For a $N \times N$ Hermitian or symmetric matrix $M$, 
for any $z \in \C\setminus {\rm Spect}(M)$, 
we denote by $G(z) := (zI_N-M)^{-1}$ the resolvent of $M$.\\
Let $z \in \C\setminus \R$, 
\begin{itemize}
\item[(i)] $\Vert G(z) \Vert \leq |\Im z|^{-1}$ where $\Vert . \Vert$ denotes the operator norm. 
\item[(ii)] $\vert G(z)_{ij} \vert \leq |\Im z|^{-1}$ for all $i,j = 1, \ldots , N$. 
\item[(iii)] For $p\geq 2$, 
\begin{equation}\label{Gij}
\frac{1}{N} \sum_{i,j = 1}^N \vert G(z)_{ij} \vert ^p \leq (|\Im z|^{-1})^p. 
\end{equation}
\item[(iv)] The derivative with respect to $M$ of the resolvent $G(z)$ satisfies: 
$$G'_M(z).B = G(z)BG(z) ~\mbox{ for any matrix $B$}.$$
\item[(v)] Let $z \in \C$ such that $|z| > \Vert M \Vert$; we have
$$\Vert G(z) \Vert \leq \frac{1}{|z| - \Vert M \Vert}.$$
\end{itemize}
\end{lemme}
\item[-] $\tilde g_N$ denotes the Stieltjes transform 
of the probability measure $\mu_\sigma \boxplus \mu_{A_N}$. 
\item[-] When we state that some quantity $\Delta _N(z)$, $z \in \C\setminus \R$, 
is $O(\frac{1}{N^p})$, this means precisely that: 
$$\vert \Delta _N(z)\vert \leq \frac{P( | \Im z |^{-1}) }{N^p},$$
for some polynomial $P$ with nonnegative coefficients which is independent of $N$. 
\item[-] For any set $S$ in $\R$, we denote the set 
$\{x \in \R , \mbox{dist}(x,S)\leq \epsilon \}$ (resp. $\{x \in \R , \mbox{dist}(x,S) < \epsilon \}$) 
by $S+[-\epsilon , +\epsilon ]$ (resp. $S+(-\epsilon , +\epsilon )$). 
\end{itemize}
\section{Free convolution}\label{freeconv}

\subsection{Definition and subordination property}

Let $\tau $ be a probability measure on $\R$. 
Its Stieltjes transform $g_\tau $ is analytic on the complex upper half-plane $\C^+$. 
There exists a domain $$D_{\alpha , \beta } = \{ u+iv \in \C, |u| < \alpha v, v > \beta \}$$
on which $g_\tau $ is univalent. 
Let $K_\tau $ be its inverse function, defined on $g_\tau (D_{\alpha , \beta })$, and 
$$R_\tau (z) = K_\tau (z) - \frac{1}{z}.$$
Given two probability measures $\tau $ and $\nu $, 
there exists a unique probability measure $\lambda $ such that
$$R_\lambda = R_\tau + R_\nu $$
on a domain where these functions are defined. 
The probability measure $\lambda $ is called 
the free convolution of $\tau $ and $\nu $ and denoted by $\tau \boxplus \nu $. 

The free convolution of probability measures has an important property, 
called subordination, which can be stated as follows: 
let $\tau $ and $\nu $ be two probability measures on $\R$; 
there exists an analytic map $F: \C^+ \rightarrow \C^+$
such that $$\forall z \in \C^+ ,  ~~~~g_{\tau \boxplus \nu}(z)= g_\nu (F(z)).$$
This phenomenon was first observed by D. Voiculescu under a genericity assumption in \cite{Voiculescu93}, 
and then proved in generality in \cite{Biane98} Theorem 3.1. 
Later, a new proof of this result was given in \cite{BelBer07}, 
using a fixed point theorem for analytic self-maps of the upper half-plane. 

\subsection{Free convolution by a semicircular distribution}

In \cite{Biane97b}, P. Biane provides a deep study of the free convolution by a semicircular distribution. 
We first recall here some of his results that will be useful in our approach.

\noindent Let $\nu $ be a probability measure on $\R$. 
P. Biane \cite{Biane97b} introduces the set 
$$\Omega _{\sigma , \nu }:=\{ u+iv \in \C^+, v > v_{\sigma , \nu }(u)\},$$ 
where the function $v_{\sigma , \nu }: \R \rightarrow \R^+$ is defined by 
$$v_{\sigma , \nu }(u) = \inf \left\{v \geq 0, \int_{\R} \frac{d\nu (x)}{(u-x)^2+v^2} \leq \frac{1}{\sigma ^2}\right\}$$ 
and proves the following

\begin{proposition}\cite{Biane97b}\label{homeo} 
The map 
$$H_{\sigma , \nu }: z \longmapsto z+\sigma ^2 g_\nu (z)$$ 
is a homeomorphism from $\overline{\Omega _{\sigma , \nu }}$ to $\C^+ \cup \R$ 
which is conformal from $\Omega _{\sigma , \nu }$ onto $\C^+$. 
Let $F_{\sigma , \nu }: \C^+ \cup \R \rightarrow \overline{\Omega _{\sigma ,\nu }}$ 
be the inverse function of $H_{\sigma , \nu }$. 
One has, 
$$\forall z \in \C^+,  ~~~~g_{\mu _\sigma \boxplus \nu }(z)= g_{\nu }(F_{\sigma , \nu }(z))$$
and then 
\begin{equation}\label{Fulton98} 
F_{\sigma , \nu }(z)=z-\sigma ^2g_{\mu _\sigma \boxplus \nu }(z).
\end{equation}
\end{proposition}

\noindent 
Note that in particular the Stieltjes transform 
$\tilde g_N$ of $\mu _\sigma \boxplus \mu _{A_N}$ satisfies 
\begin{equation} \label{eqtilde} 
\forall z \in \C ^+,  ~~~~\tilde g_N(z) = g_{\mu _{A_N}} (z - \sigma ^2 \tilde g_N(z)).
\end{equation}

\noindent 
Considering $H_{\sigma , \nu }$ as an analytic map defined in the whole upper half-plane $\C^+$, 
it is clear that
\begin{equation}\label{connexe}\Omega _{\sigma , \nu }=H_{\sigma , \nu }^{-1}(\C^+).\end{equation}
Let us give a quick proof of  (\ref{connexe}). Let $v>0$.
Since $$\Im H_{\sigma , \nu }(u+iv)=v(1-\sigma ^2 \int_{\R} \frac{d\nu (x)}{(u-x)^2+v^2}),$$ we have
\begin{equation}\label{equivch}\Im H_{\sigma , \nu }(u+iv) > 0 \Longleftrightarrow \int_{\R} \frac{d\nu (x)}{(u-x)^2+v^2} < \frac{1}{\sigma ^2}.\end{equation}
Consequently one can easily see  that $\Omega _{\sigma , \nu }$ 
is included in $H_{\sigma , \nu }^{-1}(\C^+)$. Moreover if $u+iv \in H_{\sigma , \nu }^{-1}(\C^+)$ then (\ref{equivch}) implies that 
 $v\geq v_{\sigma , \nu }(u)$. 
If we assume that $v = v_{\sigma , \nu }(u)$, then $v_{\sigma , \nu }(u) > 0$ 
and finally $$\int_{\R} \frac{d\nu (x)}{(u-x)^2+v^2} = \frac{1}{\sigma ^2}$$
by Lemma 2 in \cite{Biane97b}. 
This is a contradiction : necessarily $v > v_{\sigma , \nu }(u)$ or, in other words, 
$u+iv \in \Omega _{\sigma , \nu }$ and we are done.

~

\noindent 
The previous results of P. Biane allow him to conclude that 
$\mu _\sigma \boxplus \nu $ is absolutely continuous with respect to the Lebesgue measure 
and to obtain the following description of the support.

\begin{theoreme} \cite{Biane97b} \label{theoBiane}
Define $\Psi _{\sigma , \nu }: \R \rightarrow \R$ by: 
$$\Psi _{\sigma , \nu }(u)=H_{\sigma , \nu }(u+iv_{\sigma , \nu }(u)) = 
u+\sigma ^2 \int_{\R} \frac{(u-x)d\nu (x)}{(u-x)^2+v_\sigma (u)^2}.$$
$\Psi _{\sigma , \nu }$ is a homeomorphism and, at the point $\Psi _{\sigma , \nu }(u)$, 
the measure $\mu _\sigma \boxplus \nu $ has a density given by
$$p_{\sigma , \nu }(\Psi _{\sigma , \nu }(u))=\frac{v_{\sigma , \nu }(u)}{\pi \sigma ^2}.$$
Define the set 
$$U_{\sigma , \nu }:= \left\{ u \in \R, \int_{\R} \frac{d\nu (x)}{(u-x)^2} > \frac{1}{\sigma ^2} \right\} = 
\left\{ u \in \R , v_{\sigma , \nu }(u) > 0 \right\}.$$
The support of the measure $\mu _\sigma \boxplus \nu $ is 
the image of the closure of the open set $U_{\sigma , \nu }$ by the homeomorphism $\Psi _{\sigma , \nu }$. 
$\Psi _{\sigma , \nu }$ is strictly increasing on $U_{\sigma , \nu }$. 
\end{theoreme}

\noindent 
Hence, 
$$\R\setminus {\rm supp}(\mu _\sigma \boxplus \nu )= 
\Psi _{\sigma , \nu } (\R\setminus \overline{U_{\sigma , \nu }}).$$
One has $\Psi _{\sigma , \nu } = H_{\sigma , \nu }$ 
on $\R\setminus \overline{U_{\sigma , \nu }}$ and 
$\Psi ^{-1}_{\sigma , \nu }=F_{\sigma , \nu }$ on $\R\setminus {\rm supp} (\mu _\sigma \boxplus \nu )$. 
In particular, we have the following description of the complement of the support: 
\begin{equation}\label{ComplSupp} 
\R\setminus {\rm supp}(\mu _\sigma \boxplus \nu ) = H_{\sigma , \nu }(\R\setminus \overline{U_{\sigma , \nu }}).
\end{equation}

Let $\nu $ be a compactly supported probability measure. 
We are going to establish a characterization of the complement of the support 
of $\mu _\sigma \boxplus \nu $ involving the support of $\nu $ and $H_{\sigma , \nu }$. 
We will need the following preliminary lemma.

\begin{lemme}\label{Lemme-InclSupport}
The support of $\nu $ is included in $\overline{U_{\sigma , \nu }}$.
\end{lemme}

\noindent {\bf Proof of Lemma \ref{Lemme-InclSupport}:} 
Let $x_0$ be in $\R\setminus \overline{U_{\sigma , \nu }}$. 
Then, there is some $\epsilon > 0$ such that 
$[x_0-\epsilon , x_0+\epsilon ]\subset \R\setminus \overline{U_{\sigma , \nu }}$. 
For any integer $n \geq 1$, we define $\alpha _k = x_0-\epsilon + 2k \epsilon /n$ 
for all $0 \leq k \leq n$. 
Then, as the sets $[\alpha _k, \alpha _{k+1}]$ are trivially contained 
in $\R\setminus \overline{U_{\sigma , \nu }}$, one has that:
$$\forall u \in [\alpha _k, \alpha _{k+1}], 
\quad \frac{1}{\sigma ^2} \geq \int_{\alpha _k}^{\alpha _{k+1}} \frac{d\nu (x)}{(u-x)^2} 
\geq \frac{\nu ([\alpha _k, \alpha _{k+1}])}{(\alpha _{k+1}-\alpha _{k})^2}.$$
This readily implies that
$$\nu ([x_0-\epsilon , x_0+\epsilon ]) \leq \sum_{k=0}^{n-1} \nu ([\alpha _k, \alpha _{k+1}]) \leq 
\frac{(2\epsilon )^2}{\sigma ^2n}.$$
Letting $n \to \infty $, we get that $\nu ([x_0-\epsilon , x_0+\epsilon ])=0$, 
which implies that $x_0 \in \R\setminus {\rm supp}(\nu )$. \, $\Box$ 

~

From the  continuity and strict convexity of the function 
$u \vers \int_{\R} \frac{d\nu (x)}{(u-x)^2}$ on $\R\setminus {\rm supp}(\nu )$, it follows that
\begin{equation}\label{caract} 
\overline{U_{\sigma , \nu }} = {\rm supp}(\nu ) \cup 
\{ u \in \R\setminus {\rm supp}(\nu ), \int_{\R} \frac{d\nu (x)}{(u-x)^2} 
\geq \frac{1}{\sigma ^2} \}
\end{equation}
and
$$\R\setminus \overline{U_{\sigma , \nu }} = 
\{u \in \R\setminus {\rm supp}(\nu ), \int_{\R} \frac{d\nu (x)}{(u-x)^2} < \frac{1}{\sigma ^2}\}.$$
Now, as $H_{\sigma , \nu }$ is analytic on $\R\setminus {\rm supp}(\nu )$, 
the following characterization readily follows: 
$$\R\setminus \overline{U_{\sigma , \nu }} = 
\{ u\in \, \R\setminus {\rm supp}(\nu ), \, H_{\sigma , \nu }'(u) > 0 \}.$$
and thus, according to \eqref{ComplSupp}, we get 

\begin{proposition}\label{Caract}
$$x \in  \R\setminus {\rm supp}(\mu _\sigma \boxplus \nu ) \Leftrightarrow 
\exists u \in \R\setminus {\rm supp}(\nu ) {\rm ~such~that~} x=H_{\sigma , \nu }(u) \, , \, H_{\sigma , \nu }'(u) > 0.$$
\end{proposition}  

%
%

\begin{remarque}\label{Hcroit}
Note that $H_{\sigma , \nu }$ is strictly increasing on $\R\setminus \overline{U_{\sigma , \nu }}$ since, 
if $a < b$ are in $\R\setminus {\rm supp}(\nu )$, one has, by Cauchy-Schwarz inequality, that 
\begin{eqnarray*}
H_{\sigma , \nu }(b)-H_{\sigma , \nu }(a) & = & (b-a) \bigg [ 1 - \sigma ^2 \int_{\R}\frac{d\nu (x)}{(a-x)(b-x)} \bigg ] \\
& \geq & (b-a) \bigg [ 1 - \sigma ^2 \sqrt{ (-g'_{\nu }(a))(-g'_{\nu }(b))} \bigg ].
\end{eqnarray*}
which is nonnegative if $a$ and $b$ belong to $\R\setminus \overline{U_{\sigma , \nu }}$.
\end{remarque}

\begin{remarque}\label{compconnexes}
Each connected component of $\overline{U_{\sigma , \nu }}$ 
contains at least one connected component of ${\rm supp}(\nu )$.
\end{remarque}

Indeed, let $[s_l,t_l]$ be a connected component of $\overline{U_{\sigma , \nu }}$. 
If $s_l$ or $t_l$ is in ${\rm supp}(\nu )$, 
$[s_l,t_l]$ contains at least a connected component of ${\rm supp}(\nu)$ 
since ${\rm supp}(\nu )$ is included in $\overline{U_{\sigma ,\nu }}$. 
Now, if neither $s_l$ nor $t_l$ is in ${\rm supp}(\nu )$, 
according to (\ref{caract}), we have 
$$\int_{\R } \frac{d\nu (x)}{(s_l-x)^2} = \int_{\R } \frac{d\nu (x)}{(t_l-x)^2} =  \frac{1}{\sigma ^2}.$$ 
Assume that $[s_l,t_l]\subset \R\setminus {\rm supp}(\nu )$, 
then, by strict convexity of the function $u\longmapsto \int _\R \frac{d\nu (x)}{(u-x)^2}$ 
on $\R\setminus {\rm supp}(\nu )$, one obtains that, for any $u \in ]s_l,t_l[$, 
$$\int _\R \frac{d\nu (x)}{(u-x)^2 } < \frac{1}{\sigma ^2},$$
which leads to a contradiction. $\Box$

\begin{remarque}\label{compconnexesfini}
One can readily see that $$\overline{U_{\sigma , \nu }} \subset \{ u, {\rm dist}(u,{\rm supp}(\nu ))\leq \sigma \}$$
and deduce, since ${\rm supp}(\nu )$ is compact, that $U_{\sigma , \nu }$ is a relatively compact open set. 
Hence, $\overline{U_{\sigma , \nu }}$ has a finite number of connected components 
and may be written as the following finite disjoint union
\begin{eqnarray}{\label{DecompoU}}
\overline{U_{\sigma , \nu }} = \underset{l=m}{\overset{1}{\bigcup}} \, \bigl [s_l,t_l \bigr ] 
\quad \text{with $s_m < t_m < \ldots < s_1 < t_1$}.
\end{eqnarray}
\end{remarque}

We close this section with a proposition pointing out a relationship 
between the distribution functions of $\nu $ and $\mu _\sigma \boxplus \nu $. 

\begin{proposition} \label{palier}
Let $[s_l, t_l]$ be a connected component of $\overline{U_{\sigma , \nu }}$, then
$$(\mu _\sigma \boxplus \nu )([\Psi _{\sigma , \nu }(s_l), \Psi _{\sigma , \nu }(t_l)]) = \nu ([s_l, t_l]).$$
\end{proposition}

\noindent{\bf Proof of Proposition \ref{palier}:} 
Let $]a,b[$ be a connected component of $U_{\sigma , \nu }$. 
Since $a$ and $b$ are not atoms of $\nu $ and $\mu _\sigma \boxplus \nu $ is absolutely continuous, 
it is enough to show 
$$(\mu _\sigma \boxplus \nu )([\Psi _{\sigma , \nu }(a), \Psi _{\sigma , \nu }(b)]) = \nu([a,b]).$$
From Cauchy's inversion formula, $\mu _\sigma \boxplus \nu $ has a density given by
$p_\sigma (x) = -\frac{1}{\pi} \Im (g_\nu (F_{\nu , \sigma }(x))$ and 
$$(\mu _\sigma \boxplus \nu )([\Psi _{\sigma , \nu }(a), \Psi _{\sigma , \nu }(b)]) = 
-\frac{1}{\pi}\Im \left( \int_{\Psi _{\sigma , \nu }(a)}^{\Psi _{\sigma , \nu }(b)}
g_\nu (F_{\nu , \sigma }(x))dx\right).$$
We set $z = F_{\sigma , \nu }(x)$, then $x = H_{\sigma , \nu }(z)$ and $z = u +iv_{\sigma , \nu }(u)$.
Note that $v_{\sigma ,\nu }(u) > 0$ for $u\in ]a, b[$ 
and $v_{\sigma , \nu }(a) = v_{\sigma , \nu }(b) =0$ (see \cite{Biane97b}).
Then,
\begin{eqnarray*}
\lefteqn{(\mu _\sigma \boxplus \nu )([\Psi _{\sigma , \nu }(a), \Psi _{\sigma , \nu }(b)])}\\
&=&-\frac{1}{\pi } \Im \left( \int _a^b g_\nu (u + iv_{\sigma , \nu }(u))
H'_{\sigma ,\nu }(u + iv_{\sigma, \nu }(u))(1 + iv_{\sigma , \nu }'(u))du\right) \\
&=&-\frac{1}{\pi} \Im \left( \int _a^b g_\nu (u + iv_{\sigma , \nu }(u))
(1+ \sigma ^2g'_{\nu }(u + iv_{\sigma , \nu }(u)))(1 + iv'_{\sigma , \nu }(u))du\right) \\
&=&-\frac{1}{\pi }\left( \Im \int _a^b g_\nu (u + iv_{\sigma , \nu }(u))(1 + iv'_{\sigma , \nu }(u))du 
+ \frac{\sigma ^2}{2} \Im [g_\nu ^2(u + iv_{\sigma , \nu }(u))]_a^b \right) \\
&=&-\frac{1}{\pi } \Im \int _a^b g_\nu (u + iv_{\sigma , \nu }(u))(1 + iv'_{\sigma , \nu }(u))du = 
- \frac{1}{\pi } \Im \int _\gamma g_\nu (z)dz,
\end{eqnarray*}
where 
$$\gamma = \{z = u + iv_{\sigma , \nu }(u), u\in [a,b]\}.$$

Now, we recall that, since $a$ and $b$ are points of continuity of the distribution function of $\nu $, 
$$\nu ([a,b]) = \lim_{\epsilon \rightarrow 0} - \frac{1}{\pi } \Im \left( \int _a^b g_\nu (u+i\epsilon )du\right) = 
\lim _{\epsilon \rightarrow 0} - \frac{1}{\pi } \Im \left( \int _{\gamma _\epsilon } g_\nu (z)dz \right) ,$$
where $\gamma _\epsilon = \{ z = u + i\epsilon , u \in [a,b]\}$.
Thus, it remains to prove that:
\begin{equation} \label{contour}
\lim _{\epsilon \rightarrow 0}\left( \Im \left( \int _\gamma g_\nu (z)dz\right) 
- \Im \left( \int _{\gamma _\epsilon } g_\nu (z)dz \right) \right) = 0.
\end{equation}
Let $\epsilon > 0$ such that $\epsilon < \sup _{[a,b]} v_{\sigma , \nu }(u)$. 
We introduce the contour
$$\hat \gamma _\epsilon = \{z = u + i(v_{\sigma , \nu }(u) \wedge \epsilon ), u\in [a,b]\}.$$
From the analyticity of $g_\nu $ on $\C^+$, we have 
$$\int _\gamma g_\nu (z)dz = \int _{\hat \gamma _\epsilon } g_\nu (z)dz.$$
Let $I_\epsilon = \{u\in [a,b], v_{\sigma , \nu }(u) <  \epsilon \} = \cup C_i(\epsilon )$, 
where $C_i(\epsilon )$ are the connected components of $I_\epsilon $. 
Then, $I_\epsilon \downarrow_{\epsilon \rightarrow 0} \{a,b\}$. 
For $u \in I_\epsilon $, 
$$\vert \Im g_\nu (u+i\epsilon )\vert = \epsilon \int \frac{d\nu (x)}{(u-x)^2+\epsilon ^2} \leq 
\epsilon \int \frac{d\nu (x)}{(u-x)^2+v^2_{\sigma , \nu }(u)} \leq \frac{\epsilon }{\sigma ^2}$$
and
$$\int _{I_\epsilon }\vert \Im g_\nu (u+i\epsilon )\vert du \leq \frac{\epsilon }{\sigma ^2}(b-a).$$
On the other hand, for $u \in I_\epsilon $,
$$\vert \Im g_\nu (u+iv_{\sigma ,\nu }(u))\vert = v_{\sigma , \nu } (u) \int \frac{d\nu (x)}{(u-x)^2+v_{\sigma , \nu }(u)^2} 
\leq \frac{\epsilon }{\sigma ^2}.$$
Moreover,
$$\Re g_\nu (u+iv_{\sigma , \nu }(u))v'_{\sigma , \nu }(u) = 
\frac{\Psi _{\sigma , \nu }(u)-u}{\sigma ^2}v'_{\sigma , \nu }(u)$$
and
\begin{eqnarray*}
\int _{I_\epsilon }\Re g_\nu (u+iv_{\sigma , \nu }(u))v'_{\sigma , \nu }(u)du
&=&\int _{I_\epsilon }\frac{\Psi _{\sigma , \nu }(u)-u}{\sigma ^2}v'_{\sigma , \nu }(u)du\\
&=&\frac{1}{\sigma ^2}\sum _i[(\Psi _{\sigma , \nu }(u)-u)v_{\sigma , \nu }(u)]_{C_i(\epsilon )}\\
&&-\frac{1}{\sigma ^2}\int _{I_\epsilon }(\Psi '_{\sigma , \nu }(u)-1)v_{\sigma , \nu }(u)du,
\end{eqnarray*}
by integration by parts.
Now (see \cite{Biane97b} or Theorem \ref{theoBiane}), 
$$\int _{I_\epsilon } \Psi '_{\sigma , \nu }(u)v_{\sigma , \nu }(u)du = 
\pi \sigma ^2 (\mu _\sigma \boxplus \nu )(\Psi _{\sigma , \nu }(I_\epsilon )) \vers _{\epsilon \rightarrow 0}0.$$
$$\int _{I_\epsilon }v_{\sigma , \nu }(u)du \leq \epsilon (b-a).$$
Since $\Psi _{\sigma , \nu }$ is increasing on $[a,b]$, 
$$\sum _i[\Psi _{\sigma , \nu }(u)v_{\sigma , \nu }(u)]_{C_i(\epsilon )} \leq 
\epsilon (\Psi _{\sigma , \nu }(b)-\Psi _{\sigma , \nu }(a))$$ 
and
$$\sum _i[uv_{\sigma , \nu }(u)]_{C_i(\epsilon )} \leq \epsilon (b-a).$$
The above inequalities imply \eqref{contour}. $\Box$

\section{Approximate subordination equation for $g_N(z)$} \label{sec2}

We look for an approximative equation for $g_N(z)$ of the form \eqref{eqtilde}. 
To estimate $g_N(z)$, we first handle the simplest case 
where $W_N$ is a GUE matrix and then see how the equation is modified in the general Wigner case. 
We shall rely on an integration by parts formula.
The first integration by parts formula concerns the Gaussian case; 
the distribution $\mu $ associated to $W_N$ is a centered Gaussian distribution with variance $\sigma ^2$ 
and the resulting distribution of $X_N=W_N/\sqrt{N}$ is denoted by GUE($N,\sigma ^2/N$). 
Then, the integration by parts formula can be expressed in a matricial form.

\begin{lemme} \label{lemIPP}
Let $\Phi $ be a complex-valued $\mathcal C^1$ function on 
$(M_N(\C)_{sa})$ and $X_N \sim $ GUE($N, \frac{\sigma ^2}{N}$). 
Then,
\begin{equation} \label{IPP1}
\E[\Phi '(X_N).H] = \frac{N}{\sigma ^2}\E[\Phi (X_N)\Tr (X_N H)],
\end{equation}
for any Hermitian matrix $H$, or by linearity for $H = E_{jk}$, $1\leq j, k \leq N$, 
where $(E_{jk})_{1\leq j, k\leq N}$ is the canonical basis of the complex space of $N \times N$ matrices.
\end{lemme}

For a general distribution $\mu $, we shall use an ``approximative" integration by parts formula, 
applied to the variable $\xi = \sqrt{2}\Re ((X_N)_{kl})$ or $\sqrt{2}\Im ((X_N)_{kl})$, $k<l$, or $(X_N)_{kk}$. 
Note that for $k<l$ the derivative of $\Phi (X_N)$ with respect to $\sqrt{2}\Re ((X_N)_{kl})$ (resp. $\sqrt{2}\Im ((X_N)_{kl})$)
 is $\Phi '(X_N).e_{kl}$ (resp. $\Phi '(X_N).f_{kl}$), 
where $e_{kl} = \frac{1}{\sqrt{2}}(E_{kl} + E_{lk})$ (resp. $f_{kl} = \frac{i}{\sqrt{2}}(E_{kl} - E_{lk})$) and for any $k$, the derivative of $\Phi (X_N)$ with respect to $(X_N)_{kk}$ is  $\Phi '(X_N).E_{kk}$.

\begin{lemme} \label{lem1}
Let $\xi$ be a real-valued random variable such that  $\E(\vert \xi \vert ^{p+2}) < \infty $. 
Let $\phi $ be a function from $\R$ to $\C$ such that the first $p+1$ derivatives are continuous and bounded.
Then, 
\begin{equation}\label{IPP2}
\E (\xi \phi (\xi )) = \sum_{a=0}^p \frac{\kappa _{a+1}}{a!}\E (\phi ^{(a)}(\xi )) + \epsilon , 
\end{equation} 
where $\kappa _{a}$ are the cumulants of $\xi $, 
$\vert \epsilon \vert \leq C \sup_t \vert \phi ^{(p+1)}(t)\vert \E (\vert \xi \vert ^{p+2})$, 
$C$ only depends on $p$.
\end{lemme}

\noindent 
Let $U (=U(N))$ be a unitary matrix such that 
$$A_N= U^*{\rm diag}(\gamma _1, \ldots , \gamma _N)U$$
and let $G$ stand for $G_N(z)$. 
Consider $\tilde G = UGU^*$. 
We describe the approach in the Gaussian case and 
present the corresponding results in the general Wigner case 
but detail some technical proofs in the Appendix.
\vspace{.3cm}

\noindent
{\bf a) Gaussian case:} 
We apply \eqref{IPP1} to $\Phi (X_N) = G_{jl}$ , $H= E_{il}$, $1\leq i,j,l \leq N$, 
and then take $\frac{1}{N}\sum_l$ to obtain, 
using the resolvent equation $GX_N = -I + zG - GA_N$ (see \cite{CDF09}), 
$$Z_{ji}:= \sigma ^2 \E [G_{ji}\tr _N(G)] + \delta _{ij} - z\E (G_{ji}) + \E [(GA_N)_{ji}] = 0.$$
Now, let $1\leq k,p \leq N$ and consider the sum $\sum_{i,j}U_{ik}^*U_{pj}Z_{ji}$.
We obtain from the previous equation
\begin{equation} \label{}
\sigma ^2\E [\tilde G_{pk} \tr _N(G)] + \delta _{pk} - z\E (\tilde G_{pk}) + \gamma _k\E [\tilde G_{pk}] = 0.
\end{equation}
Hence, using Lemma \ref{varN2} in the Appendix stating that 
$$\vert \E [\tilde G_{pk} \tr_N(G)] - \E [\tilde G_{pk}] \E [\tr _N(G)] \vert = O(\frac{1}{N^2}),$$
we finally get the following estimation
\begin{equation}\label{Gtilde} 
\E (\tilde G_{pk}) = \frac{\delta _{pk}}{(z-\sigma ^2g_N(z) - \gamma _k)} + O( \frac{1}{N^2}),
\end{equation} 
where we use that $\vert \frac{1}{z-\sigma ^2g_N(z)-\gamma _i} \vert \leq |\Im z|^{-1}$, and then
\begin{eqnarray*}
g_N(z)&=&\frac{1}{N}\sum _{k=1}^N\E [\tilde G_{kk}] = 
\frac{1}{N}\sum _{i=1}^N\frac{1}{z-\sigma ^2g_N(z)-\gamma _k}+O( \frac{1}{N^2})\\
&=&\int _\R\frac{1}{z-\sigma ^2g_N(z)-x}d\mu _{A_N}(x)+O(\frac{1}{N^2})\\
&=&g_{\mu _{A_N}}(z-\sigma ^2g_N(z))+O(\frac{1}{N^2}).
\end{eqnarray*}
\noindent
In the Gaussian case, we have thus proved:

\begin{proposition} 
For $z \in \C ^+$, $g_N(z)$ satisfies:
\begin{equation} \label{}
g_N(z) = g_{\mu _{A_N}}(z-\sigma ^2g_N(z))+O\left(\frac{1}{N^2}\right).
\end{equation}
\end{proposition}

\vspace{.3cm}
\noindent
{\bf b) Non-Gaussian case: } 
In this case, the integration by parts formula gives the following generalization of \eqref{Gtilde}: 

\begin{lemme}\label{estimtildeG}
\begin{equation}
\E (\tilde G_{pk}) = \frac{\delta _{pk}}{(z-\sigma ^2g_N(z)-\gamma _k)}
+\frac{\kappa _4}{2N^2}\frac{\E[\tilde A(p,k)]}{(z-\sigma ^2g_N(z)-\gamma _k)}+O(\frac{1}{N^2}),
\end{equation}
where
\begin{eqnarray} \label{Apk}
\tilde A(p,k)&=&\sum _{i,j} U_{ik}^*U_{pj}\left\{ \sum _l G_{jl}G_{il}^3+\sum _l G_{ji}G_{il}G_{li}G_{ll}\right. \\
&&\left. +\sum _l G_{jl}G_{ii}G_{li}G_{ll}+\sum _l G_{ji}G_{ii}G_{ll}^2 \right\} \nonumber
\end{eqnarray}
and $\frac{1}{N^2}\tilde A(p,k)\leq C\frac{\vert \Im z\vert ^{-4}}{N}$.
\end{lemme}

\noindent {\bf Proof} 
Lemma \ref{estimtildeG} readily follows from \eqref{IPPtilde}, 
Lemma \ref{varN2} and \eqref{tildeA} established in the Appendix. $\Box$ 

~

\noindent 
Thus,
\begin{eqnarray*}
g_N(z)&=&\frac{1}{N}\sum_{k=1}^N\E [\tilde {G}_{kk}] = \frac{1}{N}\sum _{k=1}^N \frac{1}{z-\sigma ^2g_N(z)-\gamma _k}\\
&+&\frac{\kappa _4}{2N^3}\sum _{k=1}^N\frac{\E[\tilde A(k,k)]}{z-\sigma ^2g_N(z)-\gamma _k}+O(\frac{1}{N^2}).
\end{eqnarray*}

\noindent 
Let us show that the first three terms in $\frac{1}{N}\sum _k \E[\tilde A(k,k)]/(z-\sigma ^2g_N(z)-\gamma _k)$ 
coming from the decomposition \eqref{Apk} are bounded and thus give a $O(\frac{1}{N^2})$ contribution in $g_N(z)$. 
We denote by $G_D$ the diagonal matrix with $k$-th diagonal entry equal to $\frac{1}{z-\sigma ^2g_N(z)-\gamma _k}$.
\begin{eqnarray*}
\left\vert \sum_{i,j,k} U_{ik}^*U_{kj}\frac{1}{z-\sigma ^2g_N(z)-\gamma _k}\E [\sum _l G_{jl}G_{il}^3] \right\vert 
&=&\left\vert \E [\sum _{i,l}(U^*G_DUG)_{il}G_{il}^3] \right\vert \\
&\leq &|\Im z|^{-2}\E [\sum _{i,l}|G_{il}^3|]\\
&\leq &|\Im z|^{-5}N,
\end{eqnarray*}
using Lemma \ref{lem0}.
The second term is of the same kind. For the third term, we obtain 
$$|\sum _i(U^*G_DUG^2G^{(d)})_{ii}G_{ii}|\leq |\Im z|^{-5}N$$
where $G^{(d)}$ is the diagonal matrix with $l$-th diagonal entry equal to $G_{ll}$. 

\noindent 
It follows that 
$$g_N(z) = g_{\mu _{A_N}}(z-\sigma ^2g_N(z))+\frac{1}{N}\hat L_N(z)+O(\frac{1}{N^2}),$$
where 
\begin{equation} \label{defLN}
\hat L_N(z) = 
\frac{\kappa _4}{2N^2}\sum_{i,j,k,l}U^*_{ik}U_{kj}\frac{1}{z-\sigma ^2g_N(z)-\gamma _k}\E [G_{ji}G_{ii}G_{ll}^2].
\end{equation}
It is easy to see that $\hat L_N(z)$ is bounded by $C |\Im z|^{-5}$.

\begin{proposition}\label{ecritureLN}
$\hat L_N$ defined by \eqref{defLN} can be written as 
$$\hat L_N (z) = L_N(z) + O(\frac{1}{N}) \text{, where } L_N(z) = $$
\begin{equation} \label{LN}
\frac{\kappa _4}{2 N^2}\sum _{i,l}
[(G_{A_N}(z-\sigma ^2g_N(z)))^2]_{ii}[G_{A_N}(z-\sigma ^2g_N(z))]_{ii}([G_{A_N}(z-\sigma ^2g_N(z))]_{ll})^2.
\end{equation}
\end{proposition}

\noindent {\bf Proof of Proposition \ref{ecritureLN}:}\\
{\it Step 1:} We first show that for $1\leq a,b\leq N$,
\begin{equation} \label{ecritureGab}
\E [G_{ab}]=[G_{A_N}(z-\sigma ^2g_N(z))]_{ab}+O(\frac{1}{N}).
\end{equation}
From Lemma \ref{estimtildeG}, for any $1\leq p,k\leq N$,
$$\E [\tilde {G}_{pk}]=\frac{\delta _{pk}}{(z-\sigma ^2g_N(z)-\gamma _k)}
+\frac{\kappa _4}{2N^2}\frac{\E[\tilde A(p,k)]}{(z-\sigma ^2g_N(z)-\gamma _k)}+O(\frac{1}{N^2}).$$
Let $1\leq a,b\leq N$,
\begin{eqnarray*} 
\E [G_{ab}]&=&\sum _{p,k}U_{ap}^*\E [\tilde G_{pk}] U_{kb}\\
&=&\sum _{k}U_{ak}^*\frac{1}{(z-\sigma ^2g_N(z)-\gamma _k)}U_{kb}\\
&+&\frac{\kappa _4}{2N^2}\sum _{p,k}U_{ap}^*\frac{\E[\tilde A(p,k)]}{(z-\sigma ^2g_N(z)-\gamma _k)}U_{kb}\\
&+&O(\frac{1}{N}),
\end{eqnarray*}
since $\sum _{p,k}|U_{ap}^*U_{kb}|\leq N$. 
The first term in the right-hand side of the above equation is equal to $[G_{A_N}(z-\sigma ^2g_N(z))]_{ab}$. 
It remains to show that the term involving $\E[\tilde A(p,k)]$ is of order $\frac{1}{N}$. 
Let us consider the ``worst term" in the decomposition \eqref{Apk} of $\tilde{A}(p,k)$, namely the last one. 
\begin{eqnarray*}
\lefteqn{ \frac{1}{2N^2}\sum _{p,k,i,j,l}U_{ap}^*\frac{1}{(z-\sigma ^2g_N(z)-\gamma _k)}
U_{kb}U_{ik}^*U_{pj}\E [G_{ji}G_{ii}G_{ll}^2]}\\
&&\qquad \qquad = \frac{1}{2N^2}\E [\sum _{k,i,l}\frac{1}{(z-\sigma ^2g_N(z)-\gamma _k)}
U_{kb}U_{ik}^*G_{ai}G_{ii}G_{ll}^2]\\
&& \qquad \qquad = \frac{1}{2N^2}\E [\sum _{i,l}(U^*G_DU)_{ib}G_{ai}G_{ii}G_{ll}^2]\\
&& \qquad \qquad = \frac{1}{2N^2}\E [\sum _{l}(GG^{(d)}U^*G_DU)_{ab}G_{ll}^2]\leq \frac{1}{2N} |\Im z|^{-5}.
\end{eqnarray*}
 
\noindent
\vspace{.3cm} 
{\it Step 2:} $\hat L_N$ defined by \eqref{defLN} can be written as 
$$\frac{\kappa _4}{2N^2}\sum _{i,l}\E [(U^*G_DUG)_{ii}G_{ii}G_{ll}^2].$$
First notice the following bound (see Appendix)
\begin{equation} \label{prodesp}
\E [(U^*G_DUG)_{ii}G_{ii}G_{ll}^2]-\E [(U^*G_DUG)_{ii}]\E [G_{ii}]\E [G_{ll}]^2 = O(\frac{1}{N}).
\end{equation} 
Thus,
$$\hat L_N(z) = \frac{\kappa _4}{2N^2}\sum_{i,l}\E [(U^*G_DUG)_{ii}]\E [G_{ii}]\E [G_{ll}]^2 + O(\frac{1}{N}).$$
Now, note that $\E [(U^*G_DUG)_{ii}]=\E [(U^*G_D\tilde{G}U)_{ii}]$
and, according to Lemma~\ref{estimtildeG}, 
\begin{eqnarray*} 
\E [(U^*G_D\tilde GU)_{ii}]&=&\sum _{p,k}(U^*G_D)_{ip}\E [\tilde G_{pk}]U_{ki}\\
&=&(U^*G_D^2 U)_{ii}+\frac{\kappa _4}{2N^2}\sum _{p,k}(U^*G_D)_{ip}\E[\tilde A(p,k)](G_DU)_{ki}\\
&&+\sum _{p,k}(U^*G_D)_{ip}O_{pk}(\frac{1}{N^2})U_{ki}.
\end{eqnarray*}
 
\noindent 
Thus 

\noindent 
$\frac{\kappa _4}{2N^2}\sum_{i,l}\E [(U^*G_DUG)_{ii}]\E [G_{ii}]\E [G_{ll}]^2$
\begin{eqnarray}
&=&\frac{\kappa _4}{2N^2}\sum _{i,l}[(G_{A_N}(z-\sigma ^2g_N(z)))^2]_{ii}\E [G_{ii}]\E [G_{ll}]^2 \label{I1}\\
&&+\frac{\kappa _4^2}{4N^4}\sum_{i,l,p,k}(U^*G_D)_{ip}\E[\tilde A(p,k)](G_DU)_{ki}\E [G_{ii}]\E [G_{ll}]^2 \label{I2}\\
&&+\frac{1}{N^2}\sum_{i,l,p,k}(U^*G_D)_{ip}O_{pk}(\frac{1}{N^2})U_{ki}\E [G_{ii}]\E [G_{ll}]^2.\label{I3}
\end{eqnarray}
The last term \eqref{I3} can be rewritten as 
$$\frac{1}{N^2}\sum_{l,p,k}(U\E [G^{(d)}]U^*G_D)_{kp}O_{pk}(\frac{1}{N^2})\E[G_{ll}]^2,$$ 
so that one can easily see that it is a $O(\frac{1}{N} )$. 
  
\noindent 
The second term \eqref{I2} can be rewritten as 

$\frac{\kappa _4^2}{4N^4}\sum _{t,l,s}\E [G_{ll}]^2$
$$\times \left\{ [U^*G_DU\E [G^{(d)}]U^*G_DUG]_{ts}[G_{ts}^3+G_{tt}G_{st}G_{ss}]\right. $$
$$\hspace*{1cm}\left. + [U^*G_DU\E [G^{(d)}]U^*G_DUG]_{tt}[G_{ts}G_{st}G_{ss}+G_{tt}G_{ss}^2]\right\} ,$$
which is obviously a $O(\frac{1}{N})$. 

\noindent 
Hence, Proposition \ref{ecritureLN} follows by rewriting 
the first term \eqref{I1} using \eqref{ecritureGab}. $\Box$ 

~

\noindent
From the above computations, we can state the following : 

\begin{proposition} 
For $z \in \C ^+$, $g_N(z)$ satisfies:
\begin{equation} \label{appmatsub}
g_N(z) = g_{\mu _{A_N}}(z-\sigma ^2g_N(z))+\frac{1}{N}L_N(z)+O\left(\frac{1}{N^2}\right)
\end{equation}
where $L_N(z)$ is given by \eqref{LN}.
\end{proposition}

\section{Estimation of  $g_N-\tilde{g}_N$}\label{difference}

\begin{proposition} \label{estimdiff}
For $z \in \C^+$,
\begin{equation} \label{estimdiffeqn}
g_N(z)-\tilde g_N(z)+\frac{\tilde{E}_N(z)}{N} = O(\frac{1}{N^2}),
\end{equation}
where $\tilde{E}_N(z)$ is given by
\begin{equation}
\tilde{E}_N(z) = \{ \sigma ^2\tilde g_N'(z)-1\} \tilde{L}_N(z)
\end{equation}
$\text{with } \tilde{L}_N(z) = $
\begin{equation} \label{tildeLN}
\frac{\kappa _4}{2 N^2}\sum _{i,l}
[(G_{A_N}(z-\sigma ^2\tilde{g}_N(z)))^2]_{ii}[G_{A_N}(z-\sigma ^2\tilde{g}_N(z))]_{ii}([G_{A_N}(z-\sigma ^2\tilde{g}_N(z))]_{ll})^2.
\end{equation}
\end{proposition}

\noindent {\bf Proof of proposition \ref{estimdiff}:} 
First, we are going to prove that for $z\in \C^+$,
\begin{equation} \label{prediff}
g_N(z)-\tilde g_N(z)+\frac{E_N(z)}{N} = O(\frac{1}{N^2}),
\end{equation}
where $E_N(z)$ is given by
\begin{equation}
E_N(z) = \{ \sigma ^2\tilde g_N'(z)-1\} L_N(z).
\end{equation}
For a fixed $z\in \C^+$, one may write the subordination equation \eqref{eqtilde} : 
$$\tilde g_N(z) = g_{\mu _{A_N}}(F_{\sigma , \mu _{A_N}}(z)) = g_{\mu _{A_N}}(z - \sigma ^2\tilde g_N(z)),$$
and the approximative matricial subordination equation \eqref{appmatsub} : 
$$g_N(z) = g_{\mu _{A_N}}(z - \sigma ^2g_N(z)) + \frac{1}{N}L_N(z) +O\left(\frac{1}{N^2}\right).$$
The main idea is to simplify the difference $g_N(z)-\tilde g_N(z)$ by introducing 
a complex number $z'$ likely to satisfy 
\begin{equation} \label{fixedpoint}
F_{\sigma , \mu _{A_N}}(z')=z-\sigma ^2g_N(z).
\end{equation}
We know by Proposition \ref{homeo} that $F_{\sigma , \mu _{A_N}}$ is a homeomorphism 
from $\C^+$ to $\Omega _{\sigma , \mu _{A_N}}$ 
whose inverse $H_{\sigma , \mu _{A_N}}$ has an analytic continuation to the 
whole upper half-plane $\C^+$. 
Since $z - \sigma ^2g_N(z)\in \C^+$, $z'\in \C$ is 
well-defined by the formula : 
$$z':=H_{\sigma , \mu _{A_N}}(z - \sigma ^2g_N(z)).$$
One has 
\begin{eqnarray*}
z'-z&=&-\sigma ^2(g_N(z)-g_{\mu _{A_N}}(z - \sigma ^2g_N(z)))\\
    &=&-\sigma ^2\frac{L_N(z)}{N} + O(\frac{1}{N^2})\\
    &=&O(\frac{1}{N})
\end{eqnarray*}
There exists thus a polynomial $P$ with nonnegative coefficients 
such that $$|z'-z|\leq \frac{P(|\Im z|^{-1})}{N}.$$
On the one hand, if $$\frac{P(|\Im z|^{-1})}{N}\geq \frac{|\Im z|}{2},$$ 
or equivalently 
\begin{equation} \label{1=O(1/N)}
1\leq \frac{2|\Im z|^{-1}P(|\Im z|^{-1})}{N},
\end{equation}
it is enough to prove that 
\begin{equation} \label{O(1)}
g_N(z)-\tilde g_N(z)+\frac{E_N(z)}{N} = O(1).
\end{equation}
Indeed, if we assume that \eqref{1=O(1/N)} and \eqref{O(1)} hold, 
then there exists a polynomial $Q$ with nonnegative coefficients 
such that 
\begin{eqnarray*}
|g_N(z)-\tilde g_N(z)+\frac{E_N(z)}{N}|&\leq &Q(|\Im z|^{-1})\\
&\leq &Q(|\Im z|^{-1})\frac{2|\Im z|^{-1}P(|\Im z|^{-1})}{N}\\
&\leq &Q(|\Im z|^{-1})(\frac{2|\Im z|^{-1}P(|\Im z|^{-1})}{N})^2.
\end{eqnarray*}
Hence, $$g_N(z)-\tilde g_N(z)+\frac{E_N(z)}{N} = O(\frac{1}{N^2}).$$
To prove \eqref{O(1)}, one can notice that both $g_N(z)$ and $\tilde g_N(z)$ 
are bounded by $\frac{1}{|\Im z|}$, and that 
$$|E_N(z)|\leq \left\{ \frac{\sigma ^2}{|\Im z|^{2}}+1\right\} |L_N(z)|,$$
where $L_N(z)=O(1)$.\\
On the other hand, if $$\frac{P(|\Im z|^{-1})}{N}\leq \frac{|\Im z|}{2},$$
one has : 
$$ |\Im z'-\Im z|\leq |z'-z|\leq \frac{|\Im z|}{2}$$
which implies $\Im z'\geq \frac{\Im z}{2}$ and 
therefore $z'\in \C^+$. 
As a consequence of  (\ref{connexe}), 
$z-\sigma ^2g_N(z)\in \Omega_{\sigma , \mu _{A_N}}$ 
and \eqref{fixedpoint} is satisfied. Thus, 
$$|g_N(z)-\tilde g_N(z')-\frac{L_N(z)}{N}|\leq \frac{P(|\Im z|^{-1})}{N^2},$$
or, in other words, 
\begin{equation} \label{termone}
g_N(z)-\tilde g_N(z')-\frac{L_N(z)}{N}=O(\frac{1}{N^2}).
\end{equation}
On the other hand, 
\begin{eqnarray*}
\tilde g_N(z')-\tilde g_N(z)&=&(z-z')\int _\R \frac{d(\mu _\sigma \boxplus \mu _{A_N})(x)}{(z'-x)(z-x)}\\
                            &=&(z-z')\int _\R \frac{d(\mu _\sigma \boxplus \mu _{A_N})(x)}{(z-x)^2}\\
                            & &+(z-z')^2\int _\R \frac{d(\mu _\sigma \boxplus \mu _{A_N})(x)}{(z'-x)(z-x)^2}.
\end{eqnarray*}
Taking into account the estimation of $z'-z$ above, one has : 
$$(z-z')\int _\R \frac{d(\mu _\sigma \boxplus \mu _{A_N})(x)}{(z-x)^2}=
-\sigma ^2\tilde g_N'(z)\frac{L_N(z)}{N}+O(\frac{1}{N^2})$$
and
$$(z-z')^2\int _\R \frac{d(\mu _\sigma \boxplus \mu _{A_N})(x)}{(z'-x)(z-x)^2}=O(\frac{1}{N^2}).$$
Hence 
\begin{equation} \label{termtwo}
\tilde g_N(z')-\tilde g_N(z)+\sigma ^2\tilde g_N'(z)\frac{L_N(z)}{N}=O(\frac{1}{N^2}).
\end{equation}
(\ref{prediff}) follows from \eqref{termone} and \eqref{termtwo} since 
\begin{eqnarray*}
|g_N(z)-\tilde g_N(z)+\frac{E_N(z)}{N}|\leq &|g_N(z)-\tilde g_N(z')-\frac{L_N(z)}{N}|\\
&+|\tilde g_N(z')-\tilde g_N(z)+\sigma ^2\tilde g_N'(z)\frac{L_N(z)}{N}|.
\end{eqnarray*}
Now, since $E_N(z)=O(1)$, we can deduce from (\ref{prediff}) that $g_N(z)-\tilde g_N(z) =O(\frac{1}{N})$
and then that $E_N(z)-\tilde{E}_N(z)=O(\frac{1}{N})$. (\ref{estimdiffeqn}) readily follows. $ \quad \Box$
\begin{remarque} \label{SpectralASCV}
By combining the estimation proved above for the difference between $g_N$ 
and the Stieltjes transform of $\mu _\sigma \boxplus \mu _{A_N}$ 
with some classical arguments developed in \cite{PasLe03}, 
one can recover the almost sure convergence of the spectral distribution of $M_N$ 
to the free convolution $\mu _\sigma \boxplus \nu $.
\end{remarque}

\section{Inclusion of the spectrum of $M_N$ in a neighborhood of 
the support of $\mu _\sigma \boxplus \mu _{A_N}$}\label{inclusionN}

The purpose of this section is to prove the following Theorem \ref{inclusion}.

\begin{theoreme}\label{inclusion} 
$\forall \epsilon > 0$, 
$$\mathbb{P}(\text{ For all large N }, {\rm Spect}(M_N)\subset 
\{ x, {\rm dist}(x, {\rm supp}(\mu _\sigma \boxplus \mu _{A_N}))\leq \epsilon \} )=1.$$
\end{theoreme}

\noindent 
The proof still uses the ideas of \cite{HaaThor05} and \cite{Schultz05} 
but, since $\mu _\sigma \boxplus \mu _{A_N}$ depends on $N$, 
we need here to apply the inverse Stieltjes tranform to functions depending on $N$. 
Therefore we give the details of the proof 
to convince the reader that the approach still holds.

\begin{lemme}\label{LSt}
For any fixed large $N$, $\tilde{E}_N$ defined in Proposition \ref{estimdiff} 
is the Stieltjes transform of a compactly supported distribution $\Lambda _{N}$ on $\R$ 
whose support is included in the support of $\mu _\sigma \boxplus \mu _{A_N}$. 
\end{lemme}

The proof relies on the following characterization already used in \cite{Schultz05}. 

\begin{theoreme}\label{TS}\cite{Tillmann53}
\begin{itemize}
\item Let $\Lambda $ be a distribution on $\R$ with compact support. 
Define the Stieltjes transform of $\Lambda $, 
$ l:\C\setminus \R \rightarrow \C$ by 
$$l(z)=\Lambda \left( \frac{1}{z-x}\right) .$$
\noindent Then $l$ is analytic on $\C\setminus \R$
and has an analytic continuation to $\C\setminus {\rm supp}(\Lambda )$. 
Moreover
\begin{itemize}
\item[($c_1$)] $l(z)\rightarrow 0$ as $|z|\rightarrow \infty ,$
\item[($c_2$)] there exists a constant $C > 0$, 
an integer $n\in \N$ and a compact set $K\subset \R$ containing ${\rm supp}(\Lambda )$, 
such that for any $z\in \C\setminus \R$, 
$$|l(z)|\leq C\max \{ {\rm dist}(z,K)^{-n}, 1\} ,$$
\item[($c_3$)] for any $\phi \in \cal C^\infty (\R, \R)$ with compact support
$$\Lambda (\phi )=-\frac{1}{\pi }\lim _{y\rightarrow 0^+}\Im \int _\R\phi (x)l(x+iy)dx.$$
\end{itemize}
\item Conversely, if $K$ is a compact subset of $\R$ 
and if $l:\C \setminus K\rightarrow \C$ is an analytic function
satisfying ($c_1$) and ($c_2$) above, 
then $l$ is the Stieltjes transform of a compactly supported distribution $\Lambda $ on $\R$. 
Moreover, ${\rm supp}(\Lambda )$ is exactly the set of singular points of $l$ in $K$. 
\end{itemize}
\end{theoreme}

We use here the notations and results of Section \ref{freeconv}. 
If $u\in \R$ is not in the support of $\mu _\sigma \boxplus \mu _{A_N}$, 
according to \eqref{ComplSupp}, $u-\sigma ^2\tilde g_N(u)=F_{\sigma , \mu _{A_N}}(u)$ 
belongs to $\R\setminus \overline{U_{\sigma , \mu _{A_N}}}$ 
and then cannot belong to ${\rm Spect}(A_N)$ 
since $\mbox{Spect}(A_N)\subset U_{\sigma , \mu _{A_N}}$.
Hence the singular points of $\tilde{E}_N$ are included in the support of $\mu _\sigma \boxplus \mu _{A_N}$. 

\noindent 
Now, we are going to show that for any fixed large $N$, 
$\tilde{E}_N$ satisfies ($c_1$) and ($c_2$) of Theorem \ref{TS}. 
Let $C > 0$ be such that, for all large $N$, 
$\mbox{supp}(\mu _\sigma \boxplus \mu _{A_N})\subset [-C;C]$ and $\mbox{supp}(\mu _{A_N})\subset [-C;C].$

Let $\alpha > C+\sigma $. 
For any $z\in \C$ such that $|z| > \alpha $, 
$$\vert \sigma ^2\tilde g_N(z)\vert \leq \frac{\sigma ^2}{\vert z\vert -C}\leq \frac{\sigma ^2}{\alpha -C} 
< \frac{(\alpha -C)^2}{\alpha -C} = \alpha -C$$
and 
$$\vert z-\sigma ^2\tilde g_N(z)\vert \geq \Big \vert \vert z \vert -\vert \sigma ^2\tilde g_N(z)\vert \Big \vert 
> \vert z \vert - (\alpha -C) > C.$$
Thus we get that for any $z\in \C$ such that $|z| > \alpha $, 
\begin{eqnarray*}
\Vert G_{A_N}(z-\sigma ^2\tilde g_N(z))\Vert &\leq &\frac{1}{\vert z-\sigma ^2\tilde g_N(z)\vert -C}\\
& < &\frac{1}{\vert z\vert -(\alpha -C)-C}\\
& < &\frac{1}{\vert z\vert -\alpha }.
\end{eqnarray*}
We get readily that, for $|z| > \alpha $, 
$$\vert \tilde{E}_N(z)\vert \leq \frac{\kappa _4}{2}
\frac{1}{(|z|-\alpha )^5}\left( \frac{\sigma ^2}{(|z|-C)^2}+1 \right) .$$
Then, it is clear than $\vert \tilde{E}_N(z)\vert \rightarrow 0$ 
when $|z|\rightarrow +\infty $ and ($c_1$) is satisfied.\\
\noindent 
Now we are going to prove ($c_2$) using  the approach of \cite{Schultz05}(Lemma 5.5). 
Denote by $\mathcal {E}_N$ the convex envelope of the support of $\mu _\sigma \boxplus \mu _{A_N}$ 
and define $$K_N:=\left\{ x\in \R; {\rm dist}(x, \mathcal {E}_N)\leq 1\right\} $$
\noindent 
and $$D_N=\left\{ z\in \C; 0 < {\rm dist}(z, K_N)\leq 1\right\} .$$
\begin{itemize}
\item Let $z\in D_N\cap (\C\setminus \R)$ with $\Re (z)\in K_N$. 
We have ${\rm dist}(z, K_N)=|\Im z|\leq 1$. 
We have 
$$\vert \tilde{E}_N(z)\vert \leq \frac{\kappa _4}{2}
\left( \sigma ^2\frac{1}{\vert \Im z \vert ^2}+1 \right) \frac{1}{|\Im z|^5}.$$
Noticing that $1\leq \frac{1}{\vert \Im z\vert ^2}$, 
we easily deduce that there exists some constant $C_0$ such that 
for any $z\in D_N\cap \C\setminus \R$ with $\Re (z)\in K_N$, 
\begin{eqnarray*}
\vert \tilde{E}_N(z)\vert &\leq &C_0|\Im z|^{-7}\\
&\leq &C_0{\rm dist}(z, K_N)^{-7}\\
&\leq &C_0\max ({\rm dist}(z, K_N)^{-7}; 1).
\end{eqnarray*}
\item Let $z\in D_N\cap (\C\setminus \R)$ with $\Re (z)\notin K_N$. 
Then ${\rm dist}(z,{\rm supp}(\mu _\sigma \boxplus \mu _{A_N}))\geq 1$. 
Since $\tilde{E}_N$ is bounded on compact subsets of $\C\setminus {\rm supp}(\mu _\sigma \boxplus \mu _{A_N})$, 
we easily deduce that there exists some constant $C_1(N)$ 
such that for any $z\in D_N$ with $\Re (z)\notin K_N$, 
$$\vert \tilde{E}_N(z)\vert \leq C_1(N)\leq C_1(N)\max ({\rm dist}(z, K_N)^{-7}; 1).$$
\item Since $\vert \tilde{E}_N(z)\vert \rightarrow 0$ when $|z|\rightarrow +\infty $, 
$\tilde{E}_N$ is bounded on $\C\setminus \overline{D_N}$. 
Thus, there exists some constant $C_2(N)$ such that for any $z\in \C\setminus \overline{D_N}$, 
$$\vert \tilde{E}_N(z)\vert \leq C_2(N)=C_2(N)\max ({\rm dist}(z, K_N)^{-7}; 1).$$
\end{itemize}
Hence ($c_2$) is satisfied with $C(N)=\max (C_0, C_1(N), C_2(N))$ and $n=7$ 
and Lemma \ref{LSt} follows from Theorem \ref{TS}. $\Box$

~

\noindent {\bf Proof of Theorem \ref{inclusion}:} 
Using the inverse Stieltjes tranform, we get respectively that, 
for any $\varphi _N $ in $\cal C^\infty (\R, \R)$ with compact support, 
$$\E [\tr _N(\varphi _N(M_N))]-\int _\R \varphi _Nd(\mu _\sigma \boxplus \mu _{A_N})-\frac{\Lambda _{N}(\varphi _N)}{N}$$
$$=\frac{1}{\pi }\lim _{y\rightarrow 0^+}\Im \int _\R\varphi _N(x)r_N(x+iy)dx,$$
where $r_N(z)=\tilde g_N(z)-g_N(z)+\frac{1}{N}\tilde{E}_N(z)$ satisfies, 
according to Proposition \ref{estimdiff}, for any $z\in \C\setminus \R$, 
\begin{equation*}\label{estimgdif}
\vert r_N(z)\vert \leq \frac{1}{N^2}P(\vert \Im z \vert ^{-1}) .
\end{equation*} 
We refer the reader to the Appendix of \cite{CD07} 
where it is proved using the ideas of \cite{HaaThor05} that if  $h$ is an analytic function on $\C\setminus \R$ which satisfies
\begin{equation*}\label{nestimgdif}
\vert h(z)\vert \leq (\vert z\vert +K)^\alpha P(\vert \Im z\vert ^{-1})
\end{equation*} 
\noindent for some polynomial $P$ with nonnegative coefficients and degree $k$ and for some   numbers $K\geq 0$ and $\alpha\geq 0$, then
there exists a polynomial $Q$ such that
$$\limsup _{y\rightarrow 0^+}\vert \int _\R\varphi _N(x)h(x+iy)dx\vert $$
$$\leq \int _\R\int _0^{+\infty }\vert (1+D)^{k+1}\varphi _N(x)\vert (\vert x\vert +\sqrt{2}t+K)^\alpha Q(t)\exp(-t)dtdx$$
where $D$ stands for the derivative operator.
Hence, if there exists $K > 0$ such that, for all large $N$, 
the support of $\varphi _N$ is included in $[-K, K]$ and 
$\sup _N\sup _{x \in [-K, K]}\vert D^p\varphi _N(x)\vert =C_p < \infty$ for any $p\leq k+1$, 
dealing with $h(z) =N^2r_N(z)$, we deduce that for all large $N$,
\begin{equation*} \label{majlimsup1} 
\limsup _{y\rightarrow 0^+}\vert \int _\R \varphi _N(x)r_N(x+iy)dx\vert \leq \frac{C}{N^2}
\end{equation*} 
and then 
\begin{equation}\label{StS} 
\E [\tr _N(\varphi _N(M_N))]-\int _\R \varphi _N d(\mu _\sigma \boxplus \mu _{A_N})
-\frac{\Lambda _N(\varphi _N)}{N}=O(\frac{1}{N^2}).  
\end{equation}
Let $\rho \geq 0$ be in ${\cal C}^\infty (\R, \R)$ 
such that its support is included in $\{ \vert x\vert \leq 1\}$ and $\int \rho (x)dx=1$. 
Let $0 < \epsilon < 1$. 
Define $$\rho _{\frac{\epsilon }{2}}(x)=\frac{2}{\epsilon }\rho(\frac{2x}{\epsilon }),$$
$$K_N(\epsilon )=\{ x, {\rm dist}(x, {\rm supp}(\mu _\sigma \boxplus \mu _{A_N}))\leq \epsilon \}$$ 
and $$f_N(\epsilon )(x)=\int _\R \1 _{K_N(\epsilon )}(y)\rho _{\frac{\epsilon }{2}}(x-y)dy.$$
the function $f_{N}(\epsilon )$ is in ${\cal C}^\infty (\R, \R)$, 
$f_{N}(\epsilon )\equiv 1$ on $K_N(\frac{\epsilon }{2})$; 
its support is included in $K_N(2\epsilon )$. 
Since there exists $K$ such that, for all large $N$, the support of $\mu _\sigma \boxplus \mu _{A_N}$ 
is included in $[-K;K]$, for all large $N$ the support of $f_N(\epsilon )$ is included in $[-K-2;K+2]$ 
and for any $p > 0$, 
$$\sup _{x\in [-K-2;K+2]}\vert D^pf_N(\epsilon )(x)\vert \leq 
\sup _{x\in [-K-2;K+2]}\int _{-K-1}^{K+1} \vert D^p \rho _{\frac{\epsilon }{2}}(x-y)\vert dy \leq C_p(\epsilon ).$$
Thus, according to \eqref{StS}, 
\begin{equation} 
\E [\tr _N(f_N(\epsilon )(M_N))]-\int _\R f_N(\epsilon )d(\mu _\sigma \boxplus \mu _{A_N})
-\frac{\Lambda _N(f_N(\epsilon ))}{N}=O_{\epsilon }(\frac{1}{N^2})
\end{equation}
and 
\begin{equation}\label{prime} 
\E [\tr _N((f_N'(\epsilon ))^2(M_N))]-\int _\R (f_N'(\epsilon ))^2d(\mu _\sigma \boxplus \mu _{A_N})
-\frac{\Lambda _N((f_N'(\epsilon ))^2)}{N}=O_{\epsilon }(\frac{1}{N^2}).
\end{equation}
Moreover, following the proof of Lemma 5.6 in \cite{Schultz05}, 
one can show that $\Lambda _N(1)=0$. 
Then, the function $\psi _N(\epsilon )\equiv 1-f_N(\epsilon )$ also satisfies
\begin{equation} 
\E [\tr _N(\psi _N(\epsilon )(M_N))]-\int _\R \psi _N(\epsilon )d(\mu _\sigma \boxplus \mu _{A_N})
-\frac{\Lambda _N(\psi _N(\epsilon ))}{N}=O_{\epsilon }(\frac{1}{N^2}). 
\end{equation}
Moreover, since $\psi _N'(\epsilon )=-f_N'(\epsilon )$, 
it comes readily from \eqref{prime} that 
$$\E [\tr _N((\psi _N'(\epsilon ))^2(M_N))]-\int _\R (\psi _N'(\epsilon ))^2d(\mu _\sigma \boxplus \mu _{A_N})
-\frac{\Lambda _N((\psi _N'(\epsilon ))^2)}{N}=O_{\epsilon }(\frac{1}{N^2}).$$
Now, since $\psi _N(\epsilon )\equiv 0$ on the support of $\mu _\sigma \boxplus \mu _{A_N}$, 
we deduce that 
\begin{equation}\label{psi} 
\E [\tr _N(\psi _N(\epsilon )(M_N))]=O_{\epsilon }(\frac{1}{N^2})
\end{equation}
and 
\begin{equation} \label{psiprime}
\E [\tr _N((\psi _N'(\epsilon ))^2(M_N))]=O_{\epsilon }(\frac{1}{N^2}).
\end{equation}
By Lemma \ref{variance} (sticking to the proof of Proposition 4.7 in \cite{HaaThor05}), 
we have
$$\mathbf{V}{\left[ \tr _N(\psi _N(\epsilon )(M_N))\right] }
\leq \frac{C_\epsilon }{N^2}\E \left[ \tr _N\{ (\psi _N'(\epsilon )(M_N))^2\} \right] .$$
Hence, using \eqref{psiprime}, one can deduce that 
\begin{equation} \label{variancepsi}
\mathbf{V}{\left[ \tr _N(\psi _N(\epsilon )(M_N))\right] }=O_{\epsilon }(\frac{1}{N^4}).
\end{equation}
Set $$Z_{N, \epsilon }:=\tr _N(\psi _N(\epsilon )(M_N))$$ 
and $$\Omega _{N, \epsilon }=\{ Z_{N, \epsilon } > N^{-\frac{4}{3}}\}.$$
From \eqref{psi} and \eqref{variancepsi}, we deduce that 
$$\E \{ \vert Z_{N, \epsilon }\vert ^2\} =O_{\epsilon }(\frac{1}{N^4}).$$
Hence $$P(\Omega _{N, \epsilon })\leq N^{\frac{8}{3}}\E \{ \vert Z_{N, \epsilon }\vert ^2\}
=O_{\epsilon }(\frac{1}{N^\frac{4}{3}}).$$
By Borel-Cantelli lemma, we deduce that, almost surely for all large $N$, $Z_{N, \epsilon }\leq N^{-\frac{4}{3}}$. 
Since $Z_{N, \epsilon }\geq \1 _{\R\setminus K_N({2\epsilon })}$, 
it follows that, almost surely for all large $N$, 
the number of eigenvalues of $M_N$ which are in $\R\setminus K_N({2\epsilon })$ 
is lower than $N^{-\frac{1}{3}}$ and thus obviously has 
to be equal to zero. The proof of Theorem \ref{inclusion} is complete.$\Box$

\section{Study of $\mu _\sigma \boxplus \mu _{A_N}$}\label{supportN}

The aim of this section is to show the following inclusion 
of the support of $\mu _\sigma \boxplus \mu _{A_N}$ 
(see Theorem {\ref{Theo-InclSupportN}} below). 
To this aim, we will use the notations and results of Section \ref{freeconv}. 
We define 
\begin{eqnarray}{\label{DefTheta}}
\Theta =\{ \theta _j, \, 1 \leq j \leq J \} \quad \text{and} \quad \Theta _{\sigma , \nu }
=\Theta \cap (\R\setminus \overline{U_{\sigma , \nu }}).
\end{eqnarray}
Furthermore, for all $\theta _j \in \Theta _{\sigma , \nu }$, we set 
\begin{eqnarray}
\rho _{\theta _j}:=H_{\sigma , \nu }(\theta _j)=\theta _j+\sigma ^2g_\nu (\theta _j)
\end{eqnarray}
which is outside the support of $\mu _\sigma \boxplus \nu $ according to \eqref{ComplSupp}, 
and we define
\begin{eqnarray}
K_{\sigma , \nu }(\theta _1, \ldots , \theta _J):=
{\rm supp}(\mu _\sigma \boxplus \nu )\bigcup \left\{ \rho _{\theta _j}, \, \theta _j \in \, \Theta _{\sigma , \nu }\right\} .
\end{eqnarray}

\begin{theoreme}\label{Theo-InclSupportN}
For any $\epsilon > 0$, 
$${\rm supp}(\mu _\sigma \boxplus \mu _{A_N})\subset 
K_{\sigma , \nu }(\theta _1, \ldots , \theta _J) + (-\epsilon, \epsilon),$$ 
when $N$ is large enough.
\end{theoreme}

\noindent
Let us decompose $\mu _{A_N}$ as 
$$\mu _{A_N}=\hat \mu _{\beta , N}+\hat \mu _{\Theta , N},$$
\begin{eqnarray*}{\label{Decompo}}
\text{where $\hat \mu _{\beta , N}=\frac{1}{N}\sum _{j=1}^{N-r}\delta _{\beta _j(N)}$ 
\, and \, $\hat \mu _{\Theta , N}=\frac{1}{N}\sum_{j=1}^{J}k_j\delta _{\theta _j}$}. 
\end{eqnarray*}
In the following, we will denote by $D(x, \delta )$ the open disk centered on $x$ and with radius $\delta $. 
We begin with a trivial technical lemma we will need in the following. 

\begin{lemme}{\label{convergenceuniforme}}
Let $\cal K$ be a compact set included in $\R\setminus {\rm supp}(\nu )$. 
Then $g_{\hat \mu _{\beta , N}}'$ (which is well defined on $\cal K$ for large $N$) 
converges to $g_\nu '$ uniformly on $\cal K$. 
\end{lemme}

\noindent
{\bf Proof of Lemma \ref{convergenceuniforme}:}
We first prove  that for all $u \in \cal K$, 
\begin{eqnarray}{\label{conv1}}
-g_{\hat \mu _{\beta , N}}'(u)=\frac{1}{N}\sum_{j=1}^{N-r}\frac{1}{(u-\beta _j)^2}
\vers _{N\rightarrow +\infty }\int \frac{d\nu (x)}{(u-x)^2}=-g_\nu '(u). 
\end{eqnarray}
Let $\epsilon > 0$ be such that ${\rm dist}(\cal K, {\rm supp}(\nu ))\geq \epsilon $. 
For all $u\in \cal K$, let $h_u$ be a bounded continuous function defined on $\R$ 
which coincides with $f_u(x)=1/(u-x)^2$ on ${\rm supp}(\nu )+[-\frac{\epsilon }{2}, \frac{\epsilon }{2}]$. 
As $\max _{1\leq j\leq N-r} \, {\rm dist}(\beta _j(N), {\rm supp}(\nu ))$ tends to zero as $N\to \infty $, 
one can find $N_0$ such that, for all $N\geq N_0$, 
$\beta _j(N)\in {\rm supp}(\nu )+[-\frac{\epsilon }{2}, \frac{\epsilon }{2}]$ for all $1\leq j\leq N-r$. 
Since the sequence of measures $\hat \mu _{\beta , N}$ weakly converges to $\nu $, 
\eqref{conv1} follows, observing that $-g_{\hat \mu _{\beta , N}}'(u)=\int h_u(x) \, d\hat \mu _{\beta , N}(x)$ 
and $-g_\nu '(u)=\int h_u(x) d\nu (x)$. 

\noindent 
The uniform convergence follows from Montel's theorem, since $g_{\hat \mu _{\beta , N}}'$ 
and $g_\nu '$ are analytic on $D=\{ z\in \C , {\rm dist}(z,{\rm supp}(\nu )) > \frac{\epsilon }{2}\}$ 
and uniformly bounded on $D$ by $\frac{4}{\epsilon ^2}$ for $N\geq N_0$. $\Box$

\vspace{.3cm}
\noindent
We are now in position to give the proof of Theorem \ref{Theo-InclSupportN}. 
We recall that, from \eqref{ComplSupp}, 
\begin{eqnarray}{\label{ComplSupp-N}}
\R\setminus {\rm supp}(\mu _\sigma \boxplus \mu _{A_N})=
H_{\sigma , \mu _{A_N}}(\R\setminus \overline{U_{\sigma , \mu _{A_N}}}).
\end{eqnarray}
In the proofs, we will write for simplicity $U_N$, $H_N$ and $F_N$ instead of 
$U_{\sigma , \mu _{A_N}}$, $H_{\sigma , \mu _{A_N}}$ and $F_{\sigma , \mu _{A_N}}$ respectively.\\
The main step of the proof consists in observing the following inclusion of the open set $U_{\sigma , \mu _{A_N}}$.


\begin{lemme}{\label{Lemme-InclU-N}}
For any $\epsilon ' > 0$, 
\begin{eqnarray}{\label{InclusionU-N}}
U_{\sigma , \mu _{A_N}}\subset \, \{ u, \, {\rm dist}(u, \overline{U_{\sigma , \nu }}) < \epsilon '\} 
\cup \, \{ u, \, {\rm dist}(u, \Theta _{\sigma , \nu }) < \epsilon '\},
\end{eqnarray}
for all large $N$ (since the compact sets $\overline{U_{\sigma , \nu }}$ 
and $\Theta _{\sigma , \nu }$ are disjoint, 
the previous union is disjoint once $\epsilon '$ is small enough). 
\end{lemme}

\noindent {\bf Proof of Lemma \ref{Lemme-InclU-N}:} 
Define 
$$\mathcal F_{\epsilon '}=\{ u, \, {\rm dist}(u, \overline{U_{\sigma , \nu }})\geq \epsilon '\} 
\cap \{ u, \, {\rm dist}(u, \Theta _{\sigma , \nu })\geq \epsilon '\}.$$
We shall show that for all large $N$, $\mathcal F_{\epsilon '}\subset \R\setminus \overline{U_N}$.\\
Since $\max _{1\leq j\leq N-r} \, {\rm dist}(\beta _j(N), {\rm supp}(\nu ))\rightarrow 0\quad \text{when $N$ goes to infinity}$, 
there exists $N_0$ such that for all $N\geq N_0$, 
the $\beta _j(N)$'s are in ${\rm supp}(\nu )+(-\epsilon ', \epsilon ')$. 
Since ${\rm supp}(\nu )\subset \overline{U_{\sigma , \nu }}$, 
it is clear that for all $N\geq N_0$, 
$\mathcal F_{\epsilon '}$ is included in $\R\setminus {\rm Spect}A_N$. 
Moreover, one can readily observe that if $u$ satisfies 
${\rm dist}(u, {\rm supp}(\nu )+(-\epsilon ', \epsilon '))\geq \sigma $ 
and ${\rm dist}(u, \Theta )\geq \sigma $ then, for all $N\geq N_0$, 
$-g_{\mu _{A_N}}'(u)\leq\frac{1}{\sigma ^2}$. 
This implies that, for all $N \geq N_0$, the open set $U_N$ is included in the compact set 
$$\mathcal F_{\epsilon '}'=\{ u, \, {\rm dist}(u, {\rm supp}(\nu )+(-\epsilon ', \epsilon '))\leq \sigma \} 
\cup \{ u, \, {\rm dist}(u, \Theta )\leq \sigma \} .$$
Hence, it is sufficient to show that for $N$ large enough, the compact set 
$\mathcal K_{\epsilon '}:=\mathcal F_{\epsilon '}\cap \mathcal F_{\epsilon '}'$ 
is contained in $\R\setminus \overline{U_N}$.\\
As $\nu $ is compactly supported, the function $u \mapsto -g_\nu '(u)=\int _\R d\nu(x)/(u-x)^2$ 
is continuous on $\R\setminus {\rm supp}(\nu )$. 
Hence it reaches its bounds on the compact set $\mathcal K_{\epsilon '}$ 
(which is obviously included in $\R\setminus \overline{U_{\sigma , \nu }}$) 
so that there exists $\alpha > 0$ such that $-g_\nu '(u) \leq \frac{1}{\sigma ^2}-2\alpha $ 
for any $u$ in $\mathcal K_{\epsilon '}$.\\ 
According to Lemma \ref{convergenceuniforme}, 
there exists $N_0$ such that for all $N\geq N_0$ and for all $u$ in $\mathcal K_{\epsilon '}$, 
\begin{eqnarray}{\label{ineq1}}
|g_{\hat \mu _{\beta , N}}'(u)-g_\nu '(u)|\leq \frac{3\alpha }{4}.
\end{eqnarray}
\noindent 
At last, one can notice that $N_0$ may be chosen large enough so that 
\begin{eqnarray}
\forall N\geq N_0, \quad -g_{\hat \mu _{\Theta , N}}'(u)=
\frac{1}{N}\sum_{j=1}^J\frac{k_j}{(u-\theta _j)^2}\leq \frac{\alpha }{4}.
\end{eqnarray}
This is just because for all $u\in \mathcal F_{\epsilon '}$, 
one has that: $-g_{\hat \mu_{\Theta , N}}'(u)\leq \frac{r}{N\epsilon '^2}$ 
which converges uniformly on $\mathcal K'_{\epsilon '}$ to 0 as $N$ goes to infinity.\\
Combining all the preceding gives that, on $\mathcal K_{\epsilon '}$, 
the function $-g_{\mu _{A_N}}'$ is bounded from above by $\frac{1}{\sigma ^2}-\alpha $. 
This implies that $\mathcal K_{\epsilon '}$ is included in $\R\setminus \overline{U_{\sigma , \mu _{A_N}}}$ 
which is what we wanted to show. \, $\Box$ 

~

Now we shall establish the following inclusion. 

\begin{lemme}{\label{Lemme-supp-mu}}
For all $\epsilon > 0$, for all $\epsilon ' > 0$ small enough, 
\begin{eqnarray}{\label{Incl-2-HUN}}
\R\setminus \left( K_\sigma (\theta _1, \ldots , \theta _J)+[-\epsilon , \epsilon ]\right) \subset 
H_N \left( \{ u, \, {\rm dist}(u, \Theta _{\sigma , \nu }\cup \overline{U_{\sigma , \nu }}) > \epsilon '\} \right) ,
\end{eqnarray}
when $N$ is large enough. 
\end{lemme}

Combined with Lemma \ref{Lemme-InclU-N}, this result leads to Theorem \ref{Theo-InclSupportN}. 

~

\noindent {\bf Proof of Lemma \ref{Lemme-supp-mu}:}
According to \eqref{ComplSupp}, \eqref{DecompoU} and Remark \ref{Hcroit}, we have that
$$\R\setminus {\rm supp}(\mu _\sigma \boxplus \nu )=$$
$$\bigl] -\infty , H_{\sigma , \nu }(s_m)\bigr[ \, \bigcup \, 
\bigg( \, \underset{l=m}{\overset{2}{\bigcup }} \, \bigl] H_{\sigma , \nu }(t_l), H_{\sigma , \nu }(s_{l-1})\bigg[ 
\, \bigg) \, \bigcup \, \bigl] H_{\sigma , \nu }(t_1), +\infty \bigr[ $$
i.e. 
\begin{eqnarray}{\label{DecompoComplSup}}
{\rm supp}(\mu _\sigma \boxplus \nu )=\underset{l=m}{\overset{1}{\bigcup }}
\bigg[ H_{\sigma , \nu }(s_l), H_{\sigma , \nu }(t_l) \bigg] .
\end{eqnarray}
Note that there exists some finite integer $q$ such that, for $\epsilon $ small enough, 
$\R\setminus \left( K_\sigma (\theta _1, \ldots , \theta _J)+[-\epsilon , \epsilon ]\right) $
is the following disjoint union of intervals
$$]-\infty , h_0[\underset{i=1, \ldots ,q}{\bigcup }]k_i, h_i[\cup ]k_{q+1}, +\infty [,$$
where $h_i=H_{\sigma , \nu }(s_{p_i})-\epsilon $ and $k_{i+1}=H_{\sigma , \nu }(t_{p_i})+\epsilon $ 
for some $p_i$ or $h_i=H_{\sigma , \nu }(\theta _{j_i})-\epsilon $ 
and $k_{i+1}=H_{\sigma , \nu }(\theta _{j_i})+\epsilon $ for some $\theta _{j_i}$ in $\Theta _{\sigma , \nu }$.\\
For such an $\epsilon > 0$, since $H_{\sigma , \nu }$ coincides on 
$\R\setminus U_{\sigma , \nu }$ with the homeomorphism $\Psi_{\sigma,\nu}$ defined in Theorem \ref{theoBiane}, we can deduce in particular that
$H_{\sigma , \nu }$ is right-continuous (resp. left-continuous) at each $t_l$ (resp. $s_l$) for $1\leq l\leq m$, and 
$H_{\sigma , \nu }$ is continuous at each $\theta _i$ in $\Theta _{\sigma , \nu }$. Thus, 
there exists $\epsilon ' > 0$ such that: for all $1\leq l\leq m$, 
\begin{eqnarray}{\label{ineq2}}
H_{\sigma , \nu }(s_l-\epsilon ')\geq H_{\sigma , \nu }(s_l)-\frac{\epsilon }{2}\quad 
\text{and} \quad  H_{\sigma , \nu }(t_l+\epsilon ')\leq H_{\sigma , \nu }(t_l)+\frac{\epsilon }{2}
\end{eqnarray}
and for all $\theta _j$ in $\Theta _{\sigma , \nu }$, 
\begin{eqnarray}{\label{ineq2bis}}
H_{\sigma , \nu }(\theta _j-\epsilon ')\geq H_{\sigma , \nu }(\theta _j)-\frac{\epsilon }{2}\quad 
\text{and} \quad  H_{\sigma , \nu }(\theta _j+\epsilon ')\leq H_{\sigma , \nu }(\theta _j)+\frac{\epsilon }{2}.
\end{eqnarray}
Now $H_N$ being increasing on $\R\setminus \overline{U_N}$, for $N$ large enough,
the image by $H_N$ of 
$$\{ u, \, d(u, \Theta _{\sigma , \nu }) > \epsilon '\} \cap 
\{ u, \, d(u, \overline{U_{\sigma , \nu }}) > \epsilon '\}\subseteq \R\setminus \overline{U_N}$$
is the following disjoint union of intervals
$$]-\infty , h_0(N)[\underset{i=1, \ldots , q}{\bigcup}]k_i(N), h_i(N)[\cup ]k_{q+1}(N), +\infty [,$$
where $h_i(N)=H_N(s_{p_i}-\epsilon ')$ and $k_{i+1}(N)=H_N(t_{p_i}+\epsilon ')$ 
or $h_i(N)=H_N(\theta _{j_i}-\epsilon ')$ and $k_{i+1}(N)=H_N(\theta _{j_i}+\epsilon ')$.

\noindent 
One can see that it only remains to state that for all large $N$: $\forall 1 \leq l \leq m$, 
\begin{eqnarray}{\label{ineq3}}
H_N(s_l-\epsilon ')\geq H_{\sigma , \nu }(s_l)-\epsilon \quad 
\text{and} \quad H_N(t_l+\epsilon ')\leq H_{\sigma , \nu }(t_l)+\epsilon .
\end{eqnarray}
\begin{eqnarray}{\label{ineq4}}
H_N(\theta _i-\epsilon ')\geq H_{\sigma , \nu }(\theta _i)-\epsilon \quad 
\text{and} \quad H_N(\theta _i+\epsilon ')\leq H_{\sigma , \nu }(\theta _i)+\epsilon .
\end{eqnarray}
Moreover, as $\mu _{A_N}$ weakly converges to $\nu $, 
it is not hard to see that for all $1\leq l\leq m$, and all $\theta _i$ in $\Theta _{\sigma , \nu }$, 
$H_N(s_l-\epsilon ')$, $H_N(t_l+\epsilon ')$, $H_N(\theta _i-\epsilon ')$ and
$H_N(\theta _i+\epsilon ')$ converge as $N\to \infty $ to $H_{\sigma , \nu }(s_l-\epsilon ')$, 
$H_{\sigma , \nu }(t_l+\epsilon ')$, $H_{\sigma , \nu }(\theta _i-\epsilon ')$
and $H_{\sigma , \nu }(\theta _i+\epsilon ')$ respectively. 
So, there exists $N_0$ such that for all $N\geq N_0$: 
$H_N(s_l-\epsilon ')\geq H_{\sigma , \nu }(s_l-\epsilon ')-\frac{\epsilon }{2}$ 
and $H_N(t_l+\epsilon ')\leq H_{\sigma , \nu }(t_l+\epsilon ')+\frac{\epsilon }{2}$
as well as $H_N(\theta _i-\epsilon ')\geq H_{\sigma , \nu }(\theta _i-\epsilon ')-\frac{\epsilon }{2}$ 
and $H_N(\theta _i+\epsilon ')\leq H_{\sigma , \nu }(\theta _i+\epsilon ')+\frac{\epsilon }{2}$. 
We can then deduce \eqref{ineq3} and \eqref{ineq4} from \eqref{ineq2} and \eqref{ineq2bis}.  \, $\Box$ 


\section{Exact separation of eigenvalues}\label{exact}

Before stating the fundamental exact separation phenomenon 
between the spectrum of $M_N$ and the spectrum of $A_N$, 
we need a preliminary lemma (see Lemma \ref{gap} below).\\
From Section \ref{freeconv}, we readily deduce the following

\begin{proposition}\label{homeo1}
$$\R\setminus K_{\sigma , \nu }(\theta _1, \ldots , \theta _J)=\{ x \in \R, F_{\sigma , \nu }(x) \in \, 
\R\setminus \! \left\{ \overline{U_{\sigma , \nu }}\cup \Theta \right\} \}$$
and $F_{\sigma , \nu }$ is a homeomorphism 
from $\R\setminus K_{\sigma , \nu }(\theta _1, \ldots , \theta _J)$ 
onto $\R\setminus \! \left\{ \overline{U_{\sigma , \nu }}\cup \Theta \right\} $ with inverse $H_{\sigma , \nu }$.
\end{proposition}

\begin{remarque}\label{inc}:
For all $\hat \sigma < \sigma $, $\R\setminus \overline{U_{\sigma , \nu }}\subset \, 
\R\setminus \overline{U_{\hat \sigma , \nu }}$ 
so that it makes sense to consider the following composition of homeomorphism 
$$H_{\hat \sigma , \nu }\circ F_{\sigma , \nu }: \, 
\R\setminus K_{\sigma , \nu }(\theta _1, \ldots , \theta _J) \rightarrow 
H_{\hat \sigma , \nu }(\R\setminus \left\{ \overline{U_{\sigma , \nu }}\cup \Theta \right\} )\subset \, 
\R\setminus K_{\hat \sigma , \nu }(\theta _1, \ldots , \theta _J),$$
which is stricly increasing on each connected component of 
$\R\setminus K_{\sigma , \nu }(\theta _1, \ldots , \theta _J)$.
\end{remarque}

\begin{lemme}\label{gap}
Let $[a,b]$ be a compact set contained in $\R\setminus K_{\sigma , \nu }(\theta _1, \ldots , \theta _J)$. 
Then,
\begin{itemize}
\item[(i)] For all large $N$, $\displaystyle{[F_{\sigma , \nu }(a), F_{\sigma , \nu }(b)]}
\subset \, \R\setminus {\rm Spect}(A_N)$. 
\item[(ii)] For all $0 < \hat \sigma < \sigma $, the interval 
$[H_{\hat \sigma , \nu }(F_{\sigma , \nu }(a)), H_{\hat \sigma , \nu }(F_{\sigma , \nu }(b))]$ 
is contained in $\R\setminus K_{\hat \sigma , \nu }(\theta _1, \ldots , \theta _J)$ 
and $H_{\hat \sigma , \nu }(F_{\sigma , \nu }(b))-H_{\hat \sigma , \nu }(F_{\sigma , \nu }(a))\geq b-a$.
\end{itemize}
\end{lemme}

\noindent {\bf Proof of Lemma \ref{gap}:} 
For simplicity, we define $K_{\sigma , J}^\epsilon =K_\sigma (\theta _1, \ldots , \theta _J)+[-\epsilon , \epsilon ]$. 
As $[a,b]$ is a compact set, there exist $\epsilon > 0$ and  $\alpha > 0$ such that 
$$[a-\alpha , b+\alpha ]\subset \, \R\setminus K_{\sigma , J}^\epsilon \quad \text{and} 
\quad {\rm dist}([a-\alpha , b+\alpha ]; K_{\sigma , J}^\epsilon )\geq \alpha .$$
As before, we let $\tilde \mu _N=\mu _\sigma \boxplus \mu _{A_N}$. 
According to Theorem \ref{Theo-InclSupportN}, 
there exists some $N_0$ such that for all $N\geq N_0$, 
${\rm supp}(\tilde \mu_N)$ is contained in $K_{\sigma , J}^\epsilon $. 
Thus, using \eqref{ComplSupp} and since $F_N$ is continuous 
strictly increasing on $[a-\alpha , b+\alpha ]$, we have
\begin{eqnarray}{\label{incluFN}}
\forall N\geq N_0, \quad \displaystyle{[F_N(a-\alpha ), F_N(b+\alpha )]}\subset \, \R\setminus \overline{U}_N 
\subset \, \R\setminus {\rm Spect}(A_N).
\end{eqnarray}
As $F_{\sigma , \nu }$ is strictly increasing on the compact set $[a-\alpha , b+\alpha ]$ 
(${\rm supp}({\mu _\sigma \boxplus \nu })\subset K_{\sigma , J}^\epsilon $), 
one can consider $\delta > 0$ such that 
\begin{eqnarray}{\label{ineqF}}
F_{\sigma , \nu }(a-\alpha )\leq F_{\sigma , \nu }(a)-\delta \quad \text{and} 
\quad F_{\sigma , \nu }(b+\alpha )\geq F_{\sigma , \nu }(b)+\delta .
\end{eqnarray}
Now, the weak convergence of the probability measures $\tilde \mu _N$ 
to $\mu _\sigma \boxplus \nu $ will lead to the result, 
recalling from the definition of the subordination functions 
that for all $x \in [a-\alpha , b+\alpha ]$: $F_{\sigma , \nu }(x)=x-\sigma ^2g_{\mu _\sigma \boxplus \nu }(x)$ 
and $F_N(x)=x-\sigma ^2g_{\tilde \mu_N}(x)$ (at least for all $N\geq N_0$). 
Indeed, observing that for any $x$ in $[a-\alpha , b+\alpha ]$, 
the map $\, h:t\mapsto \frac{1}{x-t}$ is bounded on $K_{\sigma , J}^\epsilon $, 
one readily gets the simple convergence of $g_{\tilde \mu_N}$ to $g_{\mu _\sigma \boxplus \nu }$ 
as well as the one of the corresponding subordination functions, 
by considering a bounded continuous function which coincides with $h$ on $K_{\sigma , J}^\epsilon $. 
We then deduce that there exists $N_0'\geq N_0$ such that, for all $N\geq N_0'$, 
\begin{eqnarray}{\label{ineqFN}}
F_N(a-\alpha )\leq F_{\sigma , \nu }(a-\alpha )+\delta \quad \text{and} 
\quad F_N(b+\alpha )\geq F_{\sigma , \nu }(b+\alpha )-\delta .
\end{eqnarray}
Combining \eqref{incluFN}, \eqref{ineqF} and \eqref{ineqFN} proves that 
the inclusion of point (i) holds true for all $N\geq N_0'$. 

The first part of (ii) is obvious from Remark \ref{inc}. 
The second part mainly follows from the fact that $F_{\sigma , \nu }$ 
is strictly increasing on $\R\setminus {\rm supp}(\mu _\sigma \boxplus \nu )$. 
More precisely, if we set $a'=H_{\hat \sigma , \nu }(F_{\sigma , \nu }(a))$ 
and $b'=H_{\hat \sigma , \nu }(F_{\sigma , \nu }(b))$, then 
\begin{eqnarray*}
b'-a'&=& F_{\sigma , \nu }(b)-F_{\sigma , \nu }(a)+
\hat{\sigma }^2\bigl( g_\nu (F_{\sigma , \nu }(b))-g_\nu (F_{\sigma , \nu }(a)) \bigr)\\
     &\geq & F_{\sigma , \nu }(b)-F_{\sigma , \nu }(a)+\sigma ^2\bigl( g_\nu (F_{\sigma , \nu }(b))-
g_\nu (F_{\sigma , \nu }(a)) \bigr) \\
     &\geq & H_{\sigma , \nu }(F_{\sigma , \nu }(b))-H_{\sigma , \nu }(F_{\sigma , \nu }(a))=b-a
\end{eqnarray*}
since $F_{\sigma , \nu }(a) < F_{\sigma , \nu }(b)$ 
and then $g_\nu (F_{\sigma , \nu }(b))-g_\nu (F_{\sigma , \nu }(a)) < 0$. $\Box$ 

\vspace{.3cm}
The exact separation result involving the subordination function related to the free convolution 
of $\mu_\sigma$ and $\nu $ can now be stated. 
Let $[a,b]$ be a compact interval contained in $\R\setminus K_{\sigma,\nu} (\theta _1, \ldots , \theta _J)$. 
By Theorems \ref{inclusion} and \ref{Theo-InclSupportN}, 
almost surely for all large $N$, $[a,b]$ is outside the spectrum of $M_N$. 
Moreover, from Lemma \ref{gap} $(i)$, it corresponds an interval $I=[a',b']$ 
outside the spectrum of $A_N$ for all large $N$ i.e., 
with the convention that $\lambda _0(M_N)=\lambda _0(A_N)=+\infty $ and 
$\lambda _{N+1}(M_N)=\lambda _{N+1}(A_N)=-\infty $, 
there is $i_N\in \{ 0, \ldots , N\}$ such that 
\begin{equation}{\label{sep1}}
\lambda _{i_N+1}(A_N) < F_{\sigma , \nu }(a):=a' \quad 
\text{and} \quad \lambda _{i_N}(A_N) > F_{\sigma , \nu }(b):=b'. 
\end{equation}
The numbers $a$ and $a'$ (resp. $b$ and $b'$) are linked as follows: 
$$a=\rho _{a'}:= H_{\sigma , \nu }(a')=a'+\sigma ^2g_\nu (a'),$$
$$b=\rho _{b'}:= H_{\sigma , \nu }(b')=b'+\sigma ^2g_\nu (b').$$
We claim that $[a,b]$ splits the spectrum of $M_N$ exactly as $I$
splits the spectrum of $A_N$. In other words,

\begin{theoreme}{\label{Thmexact}} 
With $i_N$ satisfying \eqref{sep1}, one has
\begin{equation}\label{sep2}
\mathbb P[\lambda _{i_N+1}(M_N) < a \, \text{ and } \, \lambda _{i_N}(M_N) > b, \, \text{for all large $N$}]=1.\\
\end{equation}
\end{theoreme}

The proof closely follows the proof of Theorem 4.5 in \cite{CDF09} 
by introducing in a fit way the subordination functions or their inverses. 
For the reader's convenience, we rewrite the whole proof. 
The key idea is to introduce a continuum of matrices $M_N^{(k)}$ interpolating from $M_N$ to $A_N$: 
$$M_N^{(k)}:=\frac{\sigma _k}{\sigma }\frac{W_N}{\sqrt{N}}+A_N,$$
where $$\sigma _k^2=\sigma ^2(\frac{1}{1+kC_{a,b}}),$$
and $C_{a,b}$ being a positive constant which has to be chosen small enough 
to ensure that the matrices $M_N^{(k)}$ and $M_N^{(k+1)}$ are close enough to each other. 
More precisely, $C_{a,b}$ is chosen such that 
\begin{equation}{\label{condL}}
\max \Big( \sigma ^2C_{a,b}|g_{\mu _\sigma \boxplus \nu }(a)|; 
\sigma ^2C_{a,b}|g_{\mu _\sigma \boxplus \nu }(b)|; 3\sigma C_{a,b} \Big) < \frac{b-a}{4}. 
\end{equation}
In particular,
$\sigma _0=\sigma $ and $\sigma _k\rightarrow 0$ when $k$ goes to infinity.\\
We first prove that the intervals $[H_{\sigma _k, \nu }(F_{\sigma , \nu }(a)), H_{\sigma _k, \nu }(F_{\sigma , \nu }(b))]$ 
split respectively the spectrum of $M_N^{(k)}$ in exactly the same way.
Moreover, we also prove that for $k$ large enough, 
the interval $[H_{\sigma _k, \nu }(F_{\sigma , \nu }(a)), H_{\sigma _k, \nu }(F_{\sigma , \nu }(b))]$ 
splits the spectrum of $M_N^{(k)}$ as $[F_{\sigma , \nu }(a), F_{\sigma , \nu }(b)]$ 
splits the spectrum of $A_N$, this means roughly that we extend 
the first statement to $k=\infty $ and the result follows.

As in \cite{CDF09}, this proof is inspired by the work \cite{BaiSil99} and 
mainly relies on results on eigenvalues of the rescaled Wigner matrix $X_N$ 
combined with the following classical result (due to Weyl).

\begin{lemme}{(cf. Theorem 4.3.7 of \cite{HornJohn90})} \label{Weyl}
Let B and C be two $N \times N$ Hermitian matrices. 
For any pair of integers $j,k$ such that $1\leq j,k\leq N$ and $j+k\leq N+1$, 
we have
$$\lambda_{j+k-1}(B+C)\leq \lambda_{j}(B)+\lambda_{k}(C).$$
For any pair of integers $j,k$ such that $1\leq j,k\leq N$ and $j+k \geq N+1$, 
we have
$$\lambda_j(B)+\lambda_k(C)\leq \lambda_{j+k-N}(B+C).$$
\end{lemme}

\noindent {\bf Proof of Theorem \ref{Thmexact}:} 
Given $k\geq 0$, define 
$$a_k=H_{\sigma _k, \nu}(F_{\sigma , \nu }(a)) \, \text{and} \, b_k=H_{\sigma _k, \nu}(F_{\sigma , \nu }(b)).$$

\begin{remarque}
Note that in \cite{CDF09} where $\nu =\delta _0$, we considered
$a_k=z_{\sigma _k}(g_\sigma (a))$ where $g_\sigma $ denoted the Stieltjes transform of 
$\mu _\sigma $ and $z_{\sigma _k}$ the inverse of $g_{\sigma _k}$. 
Actually, when $\nu =\delta _0$, then $H_{\sigma _k, \nu }(z)=z+\sigma _k^2/z=z_{\sigma _k}(1/z)$ 
and $F_{\sigma , \nu }=1/g_\sigma $ so that $z_{\sigma _k}(g_\sigma )=H_{\sigma _k, \nu }(F_{\sigma , \nu })$. 
This very interpretation of the composition $z_{\sigma _k}\circ g_\sigma $ 
in terms of subordination function allows us to extend the result of exact separation to non-finite rank perturbations.
\end{remarque}

\noindent 
The last point of $(ii)$ in Lemma \ref{gap} yields $b_k-a_k\geq b-a$. 
Moreover 
\begin{eqnarray*}
a_{k+1}-a_k&=&(\sigma _{k+1}^2-\sigma _k^2)g_{\mu _\sigma \boxplus \nu }(a)\\
           &=&-C_{a,b}\frac{\sigma ^2}{(1+kC_{a,b})(1+(k+1)C_{a,b})}g_{\mu _\sigma \boxplus \nu }(a),
\end{eqnarray*}
so that $\vert a_{k+1}-a_k\vert \leq \sigma ^2C_{a,b}\vert g_{\mu _\sigma \boxplus \nu }(a)\vert .$ 
Similarly $\vert b_{k+1}-b_k\vert \leq \sigma ^2C_{a,b}\vert g_{\mu _\sigma \boxplus \nu }(b)\vert .$ 
Hence, we deduce from \eqref{condL} that
\begin{equation} \label{condL2}
|a_{k+1}-a_k| < \frac{b-a}{4} \quad \text{and} \quad |b_{k+1}-b_k| < \frac{b-a}{4}.
\end{equation}

Now, we shall show by induction on $k$ that, with probability $1$, for large $N$, 
the $M_N^{(k)}$ have respectively the same amount 
of eigenvalues to the left sides of the interval $[a_{k},b_{k}]$. 
For all $k \geq 0$, 
set $${\rm E}_k=\{ \text{no eigenvalues of $M_N^{(k)}$ in $[a_{k},b_{k}]$, for all large $N$}\} .$$
By Lemma \ref{gap} $(ii)$ and Theorems \ref{inclusion} and \ref{Theo-InclSupportN}, 
we know that $\mathbb P({\rm E}_k)=1$ for all $k$. In particular, 
one has for all $\omega \in {\rm E}_0$ and for all large $N$, 
\begin{equation}{\label{case0}} 
\exists j_N(\omega ) \in \{ 0, \ldots , N \} \text{ such that } 
\lambda _{j_N(\omega )+1}(M_N) < a \, \text{ and } \, \lambda _{j_N(\omega )}(M_N) > b.
\end{equation}
Extending the random variable $j_N$, by setting for instance $j_N:=-1$ 
on the complementary of ${\rm E}_0$, we want to show that for all $k$,
\begin{equation} {\label{casek}}
\mathbb P[\lambda _{j_N+1}(M_N^{(k)}) < a_k \, \text{ and } \, 
\lambda _{j_N}(M_N^{(k)}) > b_k, \, \text{ for all large $N$}]=1.
\end{equation}
We proceed by induction. By \eqref{case0}, this is true for $k=0$. 
Now, let us assume that \eqref{casek} holds true. 
Since $$M_N^{(k+1)}=M_N^{(k)}+(\frac{1}{\sqrt{1+(k+1)C_{a,b}}}-\frac{1}{\sqrt{1+kC_{a,b}}})X_N,$$
we can deduce from Lemma \ref{Weyl} that 
$$\lambda _{j_N+1}(M_N^{(k+1)})\leq \lambda _{j_N+1}(M_N^{(k)})+(-\lambda _N(X_N))C_{a,b}.$$
Since, for $N$ large enough, $0 < -\lambda_N(X_N)\leq 3\sigma $ almost surely, 
it follows using  \eqref{condL} that 
$$\lambda _{j_N+1}(M_N^{(k+1)}) < a_k + \frac{b-a}{4}:=\hat a_k \quad \text{a.s.}.$$
Similarly, one can show that
$$\lambda_{j_N}(M_N^{(k+1)}) > b_k-\frac{b-a}{4}:=\hat b_k \quad \text{a.s.}.$$
Inequalities \eqref{condL2} ensure that $$[\hat a_k, \hat b_k]\subset [a_{k+1}, b_{k+1}].$$
As $\mathbb P({\rm E}_{k+1})=1$, we deduce that, with probability $1$, 
$$\lambda_{j_N+1}(M_N^{(k+1)}) < a_{k+1} \, \text{ and } \, 
\lambda_{j_N}(M_N^{(k+1)}) > b_{k+1}, \quad \text{for all large $N$}.$$
This completes the proof by induction of \eqref{casek}.

Now, we are going to show that there exists $K$ large enough so that, 
for all $k \geq K$, there is exact separation of the eigenvalues of the matrices $A_N$ and $M_N^{(k)}$ i.e. 
\begin{equation}{\label{casekgen}} 
\mathbb P\big{[}\lambda _{i_N+1}(M_N^{(k)}) < a_k \, \text{ and } \, 
\lambda _{i_N}(M_N^{(k)}) > b_k, \quad \text{for all large $N$}\big{]}=1.
\end{equation}
There exists $\alpha > 0$ such that 
$[a-\alpha ; b+\alpha ]\subset \R\setminus K_{\sigma , \nu }(\theta _1, \ldots , \theta _J)$. 
Thus according to Lemma \ref{gap} (i) for all large $N$, 
$$[F_{\sigma , \nu }(a-\alpha ); F_{\sigma , \nu }(b+\alpha )]\subset \, \R\setminus {\rm Spect}(A_N).$$
Now, there exists $\epsilon ' > 0$ such that $F_{\sigma , \nu }(a-\alpha ) < F_{\sigma , \nu }(a) - \epsilon '$ 
and $F_{\sigma , \nu }(b+\alpha ) > F_{\sigma , \nu }(b)+\epsilon '$. 
It follows that, for all large $N$, 
\begin{equation}\label{FA}
\lambda _{i_N+1}(A_N) < F_{\sigma , \nu }(a)-\epsilon '\quad \text{and} 
\quad \lambda _{i_N}(A_N) > F_{\sigma , \nu }(b)+\epsilon '.
\end{equation}
Using Lemma \ref{Weyl}, \eqref{FA} and the fact that, almost surely, for all large $N$, 
$$0 < \max (-\lambda _N(X_N), \lambda _1(X_N)) < 3\sigma ,$$ 
we get the following inequalities. 

\noindent 
If $i_N < N$, for all large $N$, 
\begin{eqnarray*}
\lambda _{i_N+1}(M_N^{(k)})&\leq &\lambda _{i_N+1}(A_N)
+\frac{\sigma _k}{\sigma }\lambda _{1}(X_N)\\
&<&F_{\sigma , \nu }(a)-\epsilon '+\frac{\sigma _k}{\sigma }\lambda _{1}(X_N)\\
&=&a_k-\sigma _k^2g_{\mu _\sigma \boxplus \nu}(a)+\frac{\sigma _k}{\sigma }\lambda _{1}(X_N)-\epsilon '\\
&<&a_k-\sigma _k^2g_{\mu _\sigma \boxplus \nu}(a)+3\sigma _k-\epsilon '.
\end{eqnarray*}
If $i_N > 0$, for all large $N$, 
\begin{eqnarray*} 
\lambda _{i_N}(M_N^{(k)})&\geq &\lambda _{i_N}(A_N)+
\frac{\sigma _k}{\sigma }\lambda _N(X_N)\\
&>&F_{\sigma , \nu }(b)+\epsilon '+\frac{\sigma _k}{\sigma }\lambda _N(X_N)\\
&=&b_k-\sigma _k^2g_{\mu _\sigma \boxplus \nu }(b)+\frac{\sigma _k}{\sigma }\lambda _N(X_N)+\epsilon '\\
&>&b_k-\sigma _k^2g_{\mu _\sigma \boxplus \nu }(b)-3\sigma _k+\epsilon '.
\end{eqnarray*}

\noindent 
As $\sigma_k\to 0$ when $k\to +\infty $, there is $K$ large enough such that for all $k\geq K$, 
$$\max (|-\sigma _k^2g_{\mu _\sigma \boxplus \nu }(a)+3\sigma _k|, 
|-\sigma _k^2g_{\mu _\sigma \boxplus \nu }(b)-3\sigma _k|) < \epsilon '$$ 
and then, almost surely, for all $N$ large enough 
\begin{equation}\label{inf}
\lambda _{i_N+1}(M_N^{(k)}) < a_k \mbox{~~~~if~~} i_N < N, 
\end{equation}
\begin{equation}\label{sup} 
\quad \text{and}\quad \lambda _{i_N}(M_N^{(k)}) > b_k \mbox{~~~~if~~} i_N > 0.
\end{equation}
Since $\lambda _{N+1}(M_N^{(k)})=-\lambda _0(M_N^{(k)})=-\infty $, 
\eqref{inf} (resp. \eqref{sup}) is obviously satisfied if $i_N=N$ (resp. $i_N=0$). 
Thus, we have established that for any $i_N\in \{ 0, \ldots , N\} $ satisfying \eqref{sep1}, 
\eqref{casekgen} holds for all $k\geq K$ when $K$ is large enough. 
Comparing  this with \eqref{casek}, we deduce that $j_N=i_N$ almot surely and
$$\mathbb P\big{[}\lambda _{i_N+1}(M_N) < a \, \text{ and } \, 
\lambda _{i_N}(M_N) > b, \quad \text{ for all large $N$}\big{]}=1.$$
This ends the proof of Theorem \ref{Thmexact}. $\Box$ 

~

We readily deduce the following

\begin{corollaire}\label{thetat}
Let $\epsilon > 0$. Let us fix $u$ in $\Theta _{\sigma , \nu }\cup \{ t_l, l=1, \ldots , m\} $ 
(resp. in $\Theta _{\sigma , \nu }\cup \{ s_l, l=1, \ldots , m\} $). 
Let us choose $\delta > 0$ small enough so that for large $N$, 
$[u+\delta ; u+2\delta ]$ (resp. $[u-2\delta ; u-\delta ]$) 
is included in $(\R\setminus \overline{U_{\sigma , \nu }})\cap (\R\setminus {\rm Spect}(A_N))$ 
and for any $0\leq \delta '\leq 2\delta $, $H_{\sigma , \nu }(u+\delta ')-H_{\sigma , \nu }(u) < \epsilon $ 
(resp. $H_{\sigma , \nu }(u)-H_{\sigma , \nu }(u-\delta ') < \epsilon $). 
Let $i_N=i_N(u)$ be such that 
$$\lambda _{i_N+1}(A_N) < u+\delta \, \text{ and } \, \lambda _{i_N}(A_N) > u+2\delta $$
(resp. $\lambda_{i_N+1}(A_N) < u-2\delta \, \text{ and } \, \lambda _{i_N}(A_N) > u-\delta $). 
Then $$\mathbb P\big{[}\lambda_{i_N+1}(M_N) < H_{\sigma , \nu }(u)+\epsilon \, \text{ and } \, 
\lambda _{i_N}(M_N) > H_{\sigma , \nu }(u), \text{ for all large $N$}\big{]}=1.$$
(resp. $\mathbb P\big{[}\lambda _{i_N+1}(M_N) < H_{\sigma , \nu }(u) \, \text{ and } \, 
\lambda _{i_N}(M_N) > H_{\sigma , \nu }(u)-\epsilon \text{ for all large $N$}\big{]}=1.$)
\end{corollaire}

\section{Convergence of eigenvalues}

\noindent 
In the non-spiked case $\Theta =\emptyset $ i.e. $r=0$, 
the results of Theorems \ref{Theo-InclSupportN} and \ref{inclusion} read as: 
$\forall \epsilon > 0$, 
\begin{eqnarray}{\label{InclNonPerturb}}
\mathbb P[{\rm Spect}(M_N)\subset {\rm supp}(\mu _\sigma \boxplus \nu )+(-\epsilon ,\epsilon ), 
\, \text{for all $N$ large}]=1.
\end{eqnarray}
This readily leads to the following asymptotic result for the extremal eigenvalues.

\begin{proposition}{\label{ThmASCVNonSpike}}
Assume that the deformed model $M_N$ is without spike i.e. $r=0$. Let $k\geq 0$ be a fixed integer.\\
The first largest (resp. last smallest) eigenvalues $\lambda_{1+k}(M_N)$ 
(resp. $\lambda_{N-k} (M_N)$) converge almost surely to the right (resp. left) 
endpoint of the support of $\mu _\sigma \boxplus \nu $.
\end{proposition}

\noindent {\bf Proof of Proposition \ref{ThmASCVNonSpike}:} 
We here only focus on the convergence of the first largest eigenvalues since the other case is similar. 
Recalling that ${\rm supp}(\mu _\sigma \boxplus \nu )=\cup _{l=m}^1[H_{\sigma , \nu }(s_l),H_{\sigma , \nu }(t_l)]$, 
from \eqref{InclNonPerturb}, one has that, for all $\epsilon > 0$, 
$$\mathbb P[\limsup _N\lambda _1(M_N)\leq H_{\sigma , \nu }(t_1)+\epsilon ]=1.$$
But as $H_{\sigma , \nu }(t_1)$ is a boundary point of ${\rm supp}(\mu _\sigma \boxplus \nu )$, 
the number of eigenvalues of $M_N$ falling into $[H_{\sigma , \nu }(t_1)-\epsilon ,H_{\sigma , \nu }(t_1)+\epsilon ]$ 
tends almost surely to infinity as $N\to \infty $. Thus, almost surely, 
$$\liminf _N\lambda _{1+k}(M_N)\geq H_{\sigma , \nu }(t_1)-\epsilon .$$
The result then follows by letting $\epsilon \to 0$. $\Box$ 

~

In the spiked case where $r\geq 1$ ($\Theta \not=\emptyset $), 
the spectral measure $\mu _{M_N}$ still converges almost surely to $\mu _\sigma \boxplus \nu $. 
We shall study the impact of the spiked eigenvalues $\theta _i$'s on the local behavior of some eigenvalues of $M_N$.\\
In particular, we shall prove that once the largest spike $\theta _1$ is sufficiently big, 
the largest eigenvalue of $M_N$ jumps almost surely above the right endpoint $H_{\sigma , \nu }(t_1)$. 
Once $m\geq 2$, that is when ${\rm supp}(\mu _\sigma \boxplus \nu )$ has at least two connected components, 
we prove that there may also exist some jumps into the gap(s) of this support. 
This phenomenon holds for any $\theta _j \in \Theta_{\sigma , \nu }$.\\
For $\theta _j\not\in \Theta _{\sigma , \nu }$, that is if $\theta _j\in \overline{U_{\sigma , \nu }}$, 
two situations may occur. 
To explain this, let us consider the connected component $[s_{l_j}, t_{l_j}]$ 
of $\overline{U_{\sigma , \nu }}$ which contains $\theta _j$. 
If ${\rm supp}(\nu )\cap[\theta_j, t_{l_j}]=\emptyset $ 
(resp. ${\rm supp}(\nu )\cap [s_{l_j}, \theta _j]=\emptyset $) 
then the $k_j$ corresponding eigenvalues of $M_N$ converge almost surely 
to the corresponding boundary point $H_{\sigma , \nu }(t_{l_j})$
(resp. $H_{\sigma , \nu }(s_{l_j})$) of the support of $\mu _\sigma \boxplus \nu $. 
Otherwise, namely when $\theta _j$ is between two connected components 
of ${\rm supp}(\nu )$ included in $[s_{l_j}, t_{l_j}]$, 
the convergence occurs towards a point 
inside the (interior) of ${\rm supp}(\mu _\sigma \boxplus \nu )$.\\
Here is the precise formulation of our result. 
This is the additive analogue of the main result of 
\cite{BaiYao08b} on the almost sure convergence of the eigenvalues generated
by the spikes in a generalized spiked population model.
 
\begin{theoreme}{\label{ThmASCV}}
For each spiked eigenvalue $\theta _j$, 
we denote by $n_{j-1}+1, \ldots , n_{j-1}+k_j$ the descending ranks of $\theta _j$ among the eigenvalues of $A_N$.
\begin{itemize}
\item[{1)}] If $\theta_j \in \R\setminus \overline{U_{\sigma , \nu }}$ 
(i.e. $\in \Theta _{\sigma , \nu }$), the $k_j$ eigenvalues $(\lambda_{n_{j-1}+i}(M_N), \, 1 \leq i \leq k_j)$ 
converge almost surely outside the support of $\mu _\sigma \boxplus \nu $ 
towards $\rho _{\theta _j}=H_{\sigma , \nu }(\theta _j)$.
\item[\text{2)}] If $\theta_j \in \overline{U_{\sigma , \nu }}$ 
then we let $[s_{l_j}, t_{l_j}]$ (with $1\leq l_j\leq m$) be the connected component 
of $\overline{U_{\sigma , \nu }}$ which contains $\theta _j$.
\begin{itemize}
\item[{a)}] If $\theta_j$ is on the right (resp. on the left) of any connected component of ${\rm supp}(\nu )$ 
which is included in $[s_{l_j},t_{l_j}]$ then the $k_j$ eigenvalues $(\lambda_{n_{j-1}+i}(M_N)$, $1\leq i\leq k_j)$ 
converge almost surely to $H_{\sigma , \nu }(t_{l_j})$ (resp. $H_{\sigma , \nu }(s_{l_j})$) 
which is a boundary point of the support of $\mu _\sigma \boxplus \nu $. 
\item[{b)}] If $\theta_j$ is between two connected components of ${\rm supp}(\nu )$ 
which are included in $[s_{l_j},t_{l_j}]$ then 
the $k_j$ eigenvalues $(\lambda_{n_{j-1}+i}(M_N)$, $1\leq i\leq k_j)$ 
converge almost surely to the $\alpha _j$-th quantile of $\mu _\sigma \boxplus \nu $ 
(that is to $q_{\alpha _j}$ defined by $\alpha _j=(\mu _\sigma \boxplus \nu )(]-\infty , q_{\alpha _j}])$) 
where $\alpha _j$ is such that $\alpha _j=1-\lim _N\frac{n_{j-1}}{N}=\nu (]-\infty ,\theta _j])$.
\end{itemize}
\end{itemize}
\end{theoreme}

\noindent {\bf Proof of Theorem \ref{ThmASCV}:}
1) Choosing $u=\theta _j$ in Corollary \ref{thetat} gives, for any $\epsilon > 0$, 
\begin{eqnarray}\label{eqfin}
\rho _{\theta _j}-\epsilon \leq \lambda _{n_{j-1}+k_j}(M_N)\leq \cdots 
\leq \lambda _{n_{j-1}+1}(M_N)\leq \rho _{\theta _j}+\epsilon , \text{ for large $N$}
\end{eqnarray}
holds almost surely. Hence
$$\forall 1\leq i\leq k_j, \quad \lambda _{n_{j-1}+i}(M_N) \overset{a.s.}{\longrightarrow }\rho _{\theta _j}.$$

\noindent 
2) a) We only focus on the case where $\theta _j$ is on the right of any connected component 
of ${\rm supp}(\nu )$ which is included in $[s_{l_j}, t_{l_j}]$ 
since the other case may be considered with similar arguments. 
Let us consider the set $\{ \theta _{j_0} > \ldots > \theta_{j_p}\}$ of all the $\theta _i$'s 
being in $[s_{l_j}, t_{l_j}]$ and on the right of any connected component of ${\rm supp}(\nu )$ 
which is included in $[s_{l_j}, t_{l_j}]$. 
Note that we have for all large $N$, for any $0\leq h\leq p$, 
$$n_{j_h-1}+k_{j_h} =n_{j_h}$$ 
\noindent 
and $\theta _{j_0}$ is the largest eigenvalue of $A_N$ which is lower than $t_{l_j}$.
\noindent 
Let $\epsilon > 0$. Applying Corollary \ref{thetat} with $u=t_{l_j}$, we get that, almost surely, 
$$\lambda _{n_{j_0-1}+1}(M_N) < H_{\sigma , \nu }(t_{l_j})+\epsilon \text{ and } 
\lambda _{n_{j_0-1}}(M_N) > H_{\sigma , \nu }(t_{l_j})\text{ for all large $N$.}$$ 
Now, almost surely, the number of eigenvalues of $M_N$ being in 
$]H_{\sigma , \nu }(t_{l_j})-\epsilon , H_{\sigma , \nu }(t_{l_j})]$ 
should tend to infinity when $N$ goes to infinity. 
Since almost surely for all large $N$, $\lambda _{n_{j_0-1}}(M_N) > H_{\sigma , \nu }(t_{l_j})$ and 
$\lambda _{n_{j_0-1}+1}(M_N) < H_{\sigma , \nu }(t_{l_j})+\epsilon $, we should have 
$$H_{\sigma , \nu }(t_{l_j})-\epsilon \leq \lambda _{n_{j_p-1}+k_{j_p}}(M_N)\leq \ldots 
\leq \lambda _{n_{j_0-1}+1}(M_N) < H_{\sigma , \nu }(t_{l_j})+\epsilon.$$
Hence, we deduce that: $\forall 0\leq l\leq p$ and $\forall 1\leq i\leq k_{j_p}$, 
$\quad \lambda _{n_{j_p-1}+i}(M_N) \overset{a.s.}{\longrightarrow }H_{\sigma , \nu }(t_{l_j}).$ 
The result then follows since $j \in \{ j_0,\ldots,j_p\} $. 

\noindent 
b) Let $\alpha _j=1-\lim _N\frac{n_{j-1}}{N}=\nu (]-\infty , \theta _j])$. 
Denote by $Q$ (resp. $Q_N$) the distribution function of $\mu _\sigma \boxplus \nu $ 
(resp. of the spectral measure of $M_N$). 
Since $\mu _\sigma \boxplus \nu $ is absolutely continuous, 
$Q$ is continuous on $\R$ and strictly increasing on 
each interval $[\Psi _{\sigma , \nu }(s_l), \Psi _{\sigma , \nu }(t_l)], 1\leq l\leq m$.\\
From Proposition \ref{palier} and the hypothesis on $\theta_j$, 
$\alpha _j \in ]Q(\Psi _{\sigma , \nu }(s_{l_j})), Q(\Psi _{\sigma , \nu }(t_{l_j}))[$ 
and there exists a unique $q_j\in ]\Psi _{\sigma , \nu }(s_{l_j}), \Psi _{\sigma , \nu }(t_{l_j})[$ 
such that $Q(q_j)=\alpha _j$. Moreover, $Q$ is strictly increasing in a neighborhood of $q_i$. \\
Let $\epsilon > 0$. From the almost sure convergence of $\mu _{M_N}$ to $\mu _\sigma \boxplus \nu $, 
we deduce $$Q_N(q_j+\epsilon )\vers _{N\rightarrow \infty }Q(q_j+\epsilon ) > \alpha _j, \quad \text{a.s.}.$$
From the definition of $\alpha_j$, it follows that for large $N$, 
$N, N-1, \ldots , n_{j-1}+k_j, \ldots , n_{j-1}+1$ belong to the set $\{ k, \lambda _k(M_n)\leq q_j+\epsilon \}$ 
and thus, $$\limsup _{N\vers \infty }\lambda _{n_{j-1}+1}(M_N)\leq q_j+\epsilon .$$
In the same way, since $Q_N(q_j-\epsilon )\vers _{N\rightarrow \infty }Q(q_j-\epsilon ) < \alpha _j$,
$$\liminf _{N\vers \infty }\lambda _{n_{j-1}+k_j}(M_N)\geq q_j-\epsilon .$$
Thus, the $k_j$ eigenvalues $(\lambda_{n_{j-1}+i}(M_N)$, $1\leq i\leq k_j)$ 
converge almost surely to $q_j$. $\Box$

\section{Appendix}

We present in this appendix the different estimates on the variance used throughout the paper. 
They rely on the Poincar\'e hypothesis on the distribution $\mu $ 
of the entries of the Wigner matrix $W_N$. 
We assume that $\mu $ satisfies a Poincar\'e inequality, 
that is there exists a positive constant $C$ such that for any
$\cal C^1$ function $f: \R \rightarrow \C$  such that $f$ and $f'$ are in $L^2(\mu )$,
$$\mathbf{V}(f)\leq C\int \vert f'\vert ^2 d\mu ,$$
with $\mathbf{V}(f)=\E(\vert f-\E(f)\vert ^2)$.\\
We refer the reader to \cite{BobGot99} for a
characterization of such measures on $\R$. 
This inequality translates in the matricial case as follows:\\ 
For any matrix $M$, define $||M||_2=(\Tr (M^*M))^{\frac{1}{2}}$
the Hilbert-Schmidt norm. 
Let $\Psi : (M_N(\C)_{sa})\rightarrow \R^{N^2}$ 
be the canonical isomorphism which maps a Hermitian  matrix $M$ 
to the real parts and the imaginary parts of its entries $M_{ij}, i\leq j$. 

\begin{lemme}\label{variance} 
Let $M_N$ be the complex Wigner Deformed matrix introduced in Section \ref{intro}. 
For any $\cal C^1$ function $f:\R^{N^2} \rightarrow \C$ 
such that $f$ and its gradient $\nabla(f)$ are both polynomially bounded,
\begin{equation}\label{Poincare}
\mathbf{V}{\left[ f\circ \Psi (M_N)\right]}\leq \frac{C}{N}\E \{ \Vert \nabla \left[f\circ \Psi (M_N)\right] \Vert _2^2\} .
\end{equation}
\end{lemme}

From this Lemma and the properties of the resolvent $G$ (see Lemma \ref{lem0}), we obtain:
\begin{itemize}
\item $\mathbf{V}((G_N(z))_{ij})\leq \frac{C}{N}P(|\Im z|^{-1})$
\item $\mathbf{V}((G_N(z))_{ii}^2)\leq \frac{C}{N}P(|\Im z|^{-1})$
\item Let $H$ be a deterministic Hermitian matrix with norm $\Vert H\Vert $, then,
$$\mathbf{V}((HG_N(z))_{ii})\leq \frac{C}{N}\Vert H\Vert ^2P(|\Im z|^{-1})$$
\item $\mathbf{V}(\tr _N(G_N(z)))\leq \frac{C}{N^2}P(|\Im z|^{-1})$
\end{itemize}
where $P$ is a polynomial.
It follows that:
\begin{eqnarray*}
\E[(U^*G_DUG)_{ii}G_{ii}G_{ll}^2]&=&\E[(U^*G_DUG)_{ii}]\E[G_{ii}]\E[G_{ll}]^2 
+\frac{1}{N}P(|\Im z|^{-1}),
\end{eqnarray*}
proving \eqref{prodesp}.

We now prove 

\begin{lemme} \label{varN2} 
Let $z\in \C\setminus \R$. Then, 
$$\vert \E[\tilde G_{pk}\tr _N(G)]-\E[\tilde G_{pk}]\E[\tr _N(G)]\vert \leq \frac{P(|\Im z|^{-1})}{N^2}.$$
\end{lemme}

{\bf Proof:} The cumulant expansion gives
$$z\E(G_{ji})=\sigma ^2\E(\tr _N(G)G_{ji})+\delta _{ij}+
\E[(GA_N)_{ji}]+\frac{\kappa _4}{2N^2}\E[T(i,j)]+O_{ji}(\frac{1}{N^2}),$$
where 
\begin{eqnarray*}
T(i,j)&=&\frac{1}{3}\left\{ \frac{1}{\sqrt{2}}\sum _{l < i}\left(G_{jl}^{(3)}.(e_{li},e_{li},e_{li})+\sqrt{-1}G_{jl}^{(3)}.(f_{li},f_{li},f_{li})\right)\right. \\
      & &+\left. \frac{1}{\sqrt{2}}\sum _{l > i}\left(G_{jl}^{(3)}.(e_{il},e_{il},e_{il})-
\sqrt{-1}G_{jl}^{(3)}.(f_{il},f_{il},f_{il}) \right)\right.\\ &&\left.
+G_{jl}^{(3)}.(E_{ii},E_{ii},E_{ii})\right\} .
\end{eqnarray*}
Straightforward computations give that 
\begin{eqnarray*}
T(i,j)=&\sum _lG_{jl}G_{li}^3+\sum _lG_{ji}G_{il}G_{li}G_{ll}\\
&+\sum _lG_{jl}G_{ii}G_{li}G_{ll}+\sum_l G_{ji}G_{ii}G_{ll}^2-2G_{ii}^3G_{ji}.
\end{eqnarray*}
We now compute the sum $\sum U_{ik}^*U_{pj}\ldots $ to obtain:
\begin{eqnarray} \label{IPP4}
\nonumber (z-\gamma _k)\E[\tilde G_{pk}]=&\sigma ^2\E[\tr _N(G)\tilde G_{pk}]+
\delta _{pk}+\frac{\kappa _4}{2N^2}\E[\tilde A(p,k)]\\
&-\frac{\kappa _4}{N^2}\sum _{i,j}U^*_{ik}U_{pj}\E[G_{ii}^3G_{ji}]+
\sum _{i,j}U^*_{ik}U_{pj}O_{ji}(\frac{1}{N^2}), 
\end{eqnarray}
where $$\tilde A(p,k)=\sum _{i,j}U^*_{ik}U_{pj}A(i,j)$$ 
and 
\begin{eqnarray*}
A(i,j)=&\sum _l G_{jl}G_{li}^3+\sum _lG_{ji}G_{il}G_{li}G_{ll}\\
&+\sum _lG_{jl}G_{ii}G_{li}G_{ll}+\sum _lG_{ji}G_{ii}G_{ll}^2.
\end{eqnarray*}
Since $\frac{\kappa _4}{N^2}\sum _{i,j}U^*_{ik}U_{pj}G_{ii}^3G_{ji}=\frac{\kappa _4}{N^2}(UG(G^{(d)})^3U^*)_{pk}$, 
this term is obviously a $O(\frac{1}{N^2})$.\\
\noindent 
Let us verify the following bound for $\tilde A$:
\begin{equation} \label{tildeA}
|\frac{1}{N^2}\tilde A(p,k)|\leq C\frac{|\Im z|^{-4}}{N}.
\end{equation}
\noindent Such a bound for the first term in the decomposition of $A$ can be readily deduced from (\ref{Gij}).
We write the computation for the fourth term in the decomposition of $A$, the other two terms are similar:
\begin{eqnarray*}
\lefteqn{\frac{1}{N^2}\sum _{i,j,l}U^*_{ik}U_{pj}G_{ji}G_{ii}G_{ll}^2}\\
&&=\frac{1}{N^2}\sum _{l}(UGG^{(d)}U^*)_{pk}G_{ll}^2=O(\frac{1}{N}).
\end{eqnarray*}
We prove now that the last term in \eqref{IPP4} is of order $O(\frac{1}{N^2})$.
This term is a linear combination of terms of the form:
$$\frac{\kappa _6}{N^3}\sum _{i,j,l} U^*_{ik}U_{pj}\E[G_{jl}^{(5)}.(v_1, \ldots , v_5)],$$
where $v_u=E_{mn}$ with $(m,n)=(i,l)$ or $(m,n)=(l,i)$. 
The fifth derivative is a product of six $G$. If there are $G^2_{il}$ or $G_{il}G_{li}$ in the product, 
we can conclude thanks to Lemma \ref{lem0}. 
The only term without any $G_{il}$ is 
$$G_{ji}G_{ll}G_{ii}G_{ll}G_{ii}G_{ll}$$
which gives the contribution
$$\frac{1}{N^3}\sum _l(UG(G^{(d)})^2U^*)_{pk}G^3_{ll}=O(\frac{1}{N^2}).$$
The term with one $G_{il}$ (or $G_{li}$) will also give a contribution in $\frac{1}{N^2}$. 
Hence
\begin{equation} \label{IPPtilde}
(z-\gamma _k)\E[\tilde G_{pk}]=\sigma ^2\E[\tr _N(G)\tilde G_{pk}]+
\delta _{pk}+\frac{\kappa _4}{2N^2}\E[\tilde A(p,k)]+O(\frac{1}{N^2}). 
\end{equation}
\vspace{.3cm}
\noindent
We now apply \eqref{IPP1} (or its extension \eqref{IPP2}) 
to $\Phi (X_N)=G_{jl}G_{qq}$ and $H= E_{il}$ and take the sum in $l$. 
We obtain
\begin{eqnarray*}
z\E (G_{ji}G_{qq})=&\sigma ^2\E (\tr _N(G)G_{ji}G_{qq})+\frac{\sigma ^2}{N}\E[G_{qi}(G^2)_{jq}]+\E[G_{qq}\delta _{ij}]\\
&+\E[(GA_N)_{ji} G_{qq}]+\frac{\kappa _4}{2N^2}\E[T(i,j)G_{qq}]\\
&+\frac{\kappa _4}{2N^2}\E[B(i,j,q)]+O_{j,i}(\frac{1}{N^2}),
\end{eqnarray*}
where $B(i,j,q)$ stands for all the terms coming from 
the third derivative of the product $(G_{jl}G_{qq})$ except $G_{qq}G_{jl}^{(3)}$. 
Now, we consider $\frac{1}{N} \sum _q$ of the above equalities to obtain:
\begin{eqnarray*}
z\E (G_{ji}\tr _N(G))=&\sigma ^2\E (\tr _N(G)^2G_{ji}) +\frac{\sigma ^2}{N^2}\E[(G^3)_{ji}]+\E[\tr _N(G)\delta _{ij}]\\
&+E[(GA_N)_{ji}\tr _N(G)]+\frac{\kappa _4}{2N^2}\E[T(i,j)\tr _N(G)]\\
&+\frac{\kappa _4}{2N^2}\frac{1}{N}\sum _q\E[B(i,j,q)]+O_{j,i}(\frac{1}{N^2}).
\end{eqnarray*}
We now compute the sum $\sum U^*_{ik}U_{pj}\ldots $ and obtain
\begin{eqnarray*}
(z-\gamma _k)\E (\tilde G_{pk}\tr _N(G))=&\sigma ^2\E (\tr _N(G)^2\tilde G_{pk})+\frac{\sigma ^2}{N^2}\E[(UG^3U^*)_{pk}]\\
&+\E[\tr _N(G)\delta _{pk}]+\frac{\kappa _4}{2N^2}\E[\tilde A(p,k)\tr _N(G)]\\
&+\frac{\kappa _4}{2N^2}\frac{1}{N}\sum _q\E[\tilde B(p,k,q)]+O(\frac{1}{N^2}),
\end{eqnarray*}
where $$\tilde B(p,k,q)=\sum U^*_{ik}U_{pj}B(i,j,q)$$
and the terms $\frac{\kappa _4}{2N^2}\sum U^*_{ik}U_{pj}\E[(T(i,j)-A(i,j))\tr _N(G)]$ and 
$\sum U^*_{ik}U_{pj}O_{j,i}(\frac{1}{N^2})$ remain a $O(\frac{1}{N^2})$ 
by the same arguments used to handle the analogous terms in \eqref{IPP4}. \\
 
\noindent 
Now, consider the difference between the above equation and $g_N(z)\times$\eqref{IPP4}:

$$(z-\gamma _k)\E[(\tilde G_{pk}(\tr _N(G)-\E[\tr _N(G)])]=$$
$$\frac{\sigma ^2}{N^2}\E[(UG^3U^*)_{pk}]+\sigma ^2\E[\tr _N(G)(\tr _N(G)-\E[\tr _N(G)])\tilde G_{pk}]$$
$$+\frac{\kappa _4}{2N^2}\E[\tilde A(p,k)(\tr _N(G)-\E[\tr _N(G)])]$$
$$+\frac{\kappa _4}{2N^2}\frac{1}{N}\sum _q\E[\tilde B(p,k,q)]+O(\frac{1}{N^2})$$
and
$$(z-\gamma _k-\sigma ^2g_N(z))\E[\tilde G_{pk}(\tr _N(G)-\E[\tr _N(G)])]=$$
$$\sigma ^2\E[(\tr _N(G)-\E[\tr _N(G)])^2\tilde G_{pk}]+\frac{\sigma ^2}{N^2}\E[(UG^3U^*)_{pk}]$$
$$+\frac{\kappa _4}{2N^2}\E[\tilde A(p,k)(\tr _N(G)-\E[\tr _N(G)])]$$
$$+\frac{\kappa _4}{2N^2}\frac{1}{N}\sum _q\E[\tilde B(p,k,q)]+O(\frac{1}{N^2}).$$
We now prove that the right-hand side of the above equation is of order $\frac{1}{N^2}$. 
This is obvious for the second and first term (since $\mathbf{V}(\tr _N(G_N(z)))=O(\frac{1}{N^2})$). 
Now, we have seen that $$\frac{1}{N^2}\tilde A(p,k)\leq \frac{C|\Im z|^{-4}}{N}.$$
By Cauchy-Schwarz inequality, 
$$\frac{1}{N^2}\E[\tilde A(p,k)(\tr _N(G)-\E[\tr _N(G)])]=O(\frac{1}{N^2}).$$
It remains to study the last term 
$$\frac{1}{N^3}\sum _q\E[\tilde B(p,k,q)]=\frac{1}{N^3}\sum_{i,j,q}U^*_{ik}U_{pj}\E[B(i,j,q)].$$
This term contains derivatives of $G_{qq}$ of order $a$ with $a$ {\it strictly positive} 
($a=1,2,3$) applied to a $3$-tuple $(v_1,v_2,v_3)$ where $v_u=E_{il}$ or $E_{li}$ 
(with a product of the derivative of order $3-a$ of $G_{jl}$). 
Thus, the index $q$ appears in $\tilde B(p,k,q)$ under the form of a product 
$G_{qm}G_{nq}$ with $m, n \in \{ i, l\}$. 
Thus, the sum in $q$ will give $G^2_{nm}$. 
Moreover, the term in $j$ in the derivative appears as $G_{jm}$ 
with $m\in \{ i, l\}$ and we can do the sum in $j$ to obtain $(UG)_{pm}$. 
Thus, $\frac{1}{N^3}\sum _q\tilde B(p,k,q)$ can be written as $\frac{1}{N^3}\sum_{i,l}$ of terms of the form
$$U^*_{ik}(G^2)_{i_1j_1}(UG)_{pj_2}G_{i_3j_3}G_{i_4j_4},$$
where $i_r, j_r \in \{ i, l\}$ and $j_2=l$ for $a=3$ (no derivative in $G_{jl}$), $j_4=l$ for $a < 3$. 
As in the previous computations, either the product $G^2_{il}$ (or $G_{il}G_{li}$) appears 
and we can apply Lemma \ref{lem0} (the others terms are bounded). 
In the other cases, we can always perform one sum in $i$ (or $l$) 
and obtain  $\frac{1}{N^3}\sum _{l(\mbox{ or }i)}$ of bounded terms. 
Let us just give an example of terms which can be obtained (for $a=1$):
$$U^*_{ik}(G^2)_{li}(UG)_{pl}G_{ii}G_{ll}.$$
Then, 
$$\frac{1}{N^3}\sum_{i,l}U^*_{ik}(G^2)_{li}(UG)_{pl}G_{ii}G_{ll}=
\frac{1}{N^3}\sum_{i}U^*_{ik}(UG G^{(d)}G^2)_{pi} G_{ii}.$$
Therefore, $\frac{1}{N^3}\sum_q \E[\tilde B(p,k,q)]$ is of order $\frac{1}{N^2}$. 
This proves Lemma \ref{varN2} since 
$\vert \frac{1}{z-\gamma_k-\sigma ^2g_N(z)}\vert \leq |\Im z|^{-1}$. $\Box$
\def\cprime{$'$}

\end{document}